\documentclass[11pt,a4paper]{article}
 \usepackage[cp1251]{inputenc}
 \usepackage[english]{babel}
 \usepackage{hhline}
 \usepackage{longtable}
 \usepackage{amssymb,amsmath,amsfonts,latexsym,graphicx,epstopdf}

\usepackage[T2A]{fontenc}

 \newcommand{\be}{\begin{equation}}
\newcommand{\ee}{\end{equation}}
\newtheorem{theo}{Theorem}
\newtheorem{remark}{Remark}

\newcommand{\bt}{\begin{theo}}
\newcommand{\et}{\end{theo}}

\begin{document}

\begin{center}
 {\Large \bf A complete Lie symmetry classification of a class of
 (1+2)-dimensional reaction-diffusion-convection equations  }
\medskip

{\bf Roman Cherniha,$^{a}$\footnote{\small  Corresponding author. E-mail: r.m.cherniha@gmail.com}}
  {\bf Mykola Serov $^b$ and Yulia Prystavka $^b$}
 \\
{\it $^a$~Institute of Mathematics,  National Academy
of Sciences  of Ukraine,\\
 3, Tereshchenkivs'ka Street, Kyiv 01004, Ukraine \\
  $^b$~\it Poltava National Technical Yuri Kondratyuk University, \\24, Pershotravnevyi Prospekt, 36601 Poltava, Ukraine
}\\
 \end{center}

 \begin{abstract}

A class of nonlinear reaction-diffusion-convection equations describing various processes in physics, biology, chemistry etc. is under study in the case of  time and  two  space  variables. The group of equivalence transformations is constructed, which is applied for deriving   a  Lie symmetry classification for the class of such equations
by   the well-known algorithm.
 It is proved that the algorithm leads to  32 reaction-diffusion-convection  equations admitting nontrivial Lie symmetries.  Furthermore a set of form-preserving transformations for this class  is constructed in order to reduce this  number of the equations and obtain a complete Lie symmetry classification. As a result, the so called canonical list of all inequivalent  equations  admitting nontrivial Lie symmetry (up to any point  transformations)    and their Lie symmetries are derived. The list consists of 22 equations and it is shown  that any other reaction-diffusion-convection equation admitting a nontrivial Lie symmetry is reducible to one of these 22 equations. As a nontrivial example,
 the symmetries
derived are applied for the reduction and finding exact solutions in the case of the porous-Fisher type equation with the Burgers  term.

\end{abstract}

\section{\bf Introduction}

Nowadays, it is generally accepted   that a huge number of real processes arising in physics, biology, chemistry, material sciences, engineering, ecology, economics  etc. can be  adequately  described only  by {\it nonlinear PDEs} (or systems of such equations). The most widely used type of equations  for modeling such processes are  nonlinear reaction-diffusion-convection~(advection) equations.  In the 1970s several monographs were published, which are devoted to study and application of the nonlinear RDC equations  in physics \cite{ames-65,am-72,ko-75}, biology \cite{fife-79,mu-77} and  chemistry \cite{aris-75I,aris-75II}. In our opinion, these books had a great impact attracting many scholars  to use reaction-diffusion-convection~(RDC) for modeling real world processes  and to study their properties.
 During the last two decades many new   monographs    appeared, especially for models related to the life sciences (see,  \cite{br-03,edel-05,ku-na-ei-16, mu-02, mu-03,ok-le-01, wa-15, ch-da-2017}).

The most general class of RDC equations occurring in applications  reads as
\be\label{1-1-0} u_t=
\nabla\cdot\left(D(u)\nabla u \right)+ K(u)\cdot\nabla u +R(u). \ee
  where $u$ is the function of $t$ and $x_1,\dots,x_n$, $\nabla = \left(\frac{\partial}{\partial x_1},\dots,
\frac{\partial}{\partial x_n} \right)$, and $\cdot$ means the scalar product.
 Here   he functions $D, K$ and $R$ are related
to the  three   most common   types of  transport  mechanisms occurring in real world processes.  The diffusivity $D(u)>0$  is the main characteristic of the diffusion (heat conductivity) process, the  vector $K(u)$ typically means velocity, which can be positive and/or negative and  describes the  convective transport (in contrast to diffusion, one is not random) and the reaction term $R(u)$ describes the  kinetics process(for example,  this function presents  interaction of the population $u$ with the environment and  its birth-death rate).
Notably, any RDC equation with the constant vector $K(u)=(c_1, \dots, c_n)$ is reducible to the equation with the same structure but with $K(u)=0$
via the known Galilei transformation

\[ x'_1 = x_1 +c_1 t, , \dots, x'_n = x_n +c_n t.  \]

 The most common equations among the class of RDC equations  (\ref{1-1-0}) arising  in applications are diffusion (heat) equations. Their typical form  is
 \be\label{1-1-1} u_t=
\nabla\cdot\left(D(u)\nabla u \right). \ee
Another important case occurs if one takes into account the diffusion and the kinetics process, hence the so called reaction-diffusion (RD)  equation
\be\label{1-1-2} u_t=
\nabla\cdot\left(D(u)\nabla u \right) +R(u) \ee
is obtained.

Lie symmetries of RDC equations of the various forms were under study starting from the classical works S.Lie \cite{Lie-1881,Lie-1885},
who calculated the maximal algebra of invariance of the linear heat equation in one-dimensional space, i.e. $u_t= u_{xx}$. Much later, in the second half of 20th century, several papers were devoted to Lie symmetry classification (another common terminology is group classification) of different subclasses of  (\ref{1-1-0}) in the case $n=1$, i.e. one space variable. The most general result was derived in \cite{ch-se-98, ch-se-2006} (see also extended version in Chapter 2 of  \cite{ch-se-pl-2018}). So, the problem of Lie symmetry classification (LSC) of the class of RDC equations (\ref{1-1-0}) in the case $n=1$ has been  completely solved.

The LSC problem of (\ref{1-1-0}) is still open in the multidimensional case,  i.e. $n>1$. Lie symmetries in the simplest case, when (\ref{1-1-0}) is the multidimensional  linear heat equation  $u_t=\triangle u$, were derived in papers \cite{go-27} (case $n=2$) and \cite{ni-72} (case $n=3$). Obviously these results can be straightforwardly generalized on the arbitrary number of the space variables.

The first LSC of the  class of nonlinear heat  equations  (\ref{1-1-1}) with two space variables was derived in paper \cite{na-70} (incidentally not cited so often as \cite{do-kn-sv-83} published 13 years later). In fact, all generic extensions of Lie symmetry depending on the form of the diffusivity $D(u)$ in equation (\ref{1-1-1}) were identified in \cite{na-70}. In particular, the author has shown that the diffusivity $D(u)=u^{-1}$  leads to infinite-dimensional Lie algebra of invariance.

A complete LSC of the class of RD equations (\ref{1-1-2}) was derived \cite{do-kn-sv-83}. In particular, it was shown that the result is essentially different in cases $n=2$  and $n=3$.

Paper \cite{ed-broa-94} is devoted to LSC of the class of reaction-convection equations
\be\label{1-1-3} u_t=
\nabla\cdot\left(D(u)\nabla u \right)+ K(u)\cdot\nabla u  \ee
with the vector  $K(u)= (K^1, 0)$ and  $n=2$. All possible Lie symmetries are listed in Table 1 \cite{ed-broa-94}. It turns out that there are only five  nontrivial extensions of the principal (basic)
Lie algebra of invariance
\begin{equation}\begin{array}{l}\label{1-1-4}
A^{pr}=<X_0=\frac{\partial}{\partial t},\quad
X_1=\frac{\partial}{\partial x_1},\quad
X_2=\frac{\partial}{\partial x_2}>.
\end{array}\end{equation}
 We remind the reader that the Lie algebra (\ref{1-1-4}) is called  principal because one is the largest  common algebra of invariance for all equations of the form (\ref{1-1-3}). It should be also mentioned that  the results derived in  \cite{ed-broa-94} were generalized in paper \cite{de-iv-soph-08} on the class of equations (\ref{1-1-3}) with $n=3$, moreover the so called anisotropic case was also under study.

 To the best of our knowledge, the classical problem, which consists in  deriving
  a complete LSC of the class of RDC equations (\ref{1-1-0}) with $n>1$, is not solved at the present time. This problem cannot be solved  as a consequence or generalization of the results obtained in the papers cited above.
  Moreover, it is well-known that multidimensional (in space) case is not reducible to 1D case in Lie symmetry analysis. For example, the so called conformal exponent $k=-4/(n+2)$ in the diffusivity $D(u)=u^k$ leads to absolutely different Lie symmetries of equation (\ref{1-1-1}) for $n=1$  and $n=2$, hence the results obtained in \cite{ch-se-98, ch-se-2006} cannot be generalized on the 2D case.
Finally, it follows from the papers \cite{na-70} and  \cite{do-kn-sv-83} that the case $n=2$ is a special one comparing with $n>2$. {\it Thus, we examine here  equations of the form (\ref{1-1-0}) in the case of two space variables.}

 In Section 2,  the determining equations (DEs) for finding Lie symmetries are constructed  and the group of   equivalence transformations (ETs) is identified. It is proved that continuous ETs of the class of RDC equations (\ref{1-1-0})($n=2$) form the 10-parameter Lie group.
 In Section 3,  necessary conditions for existing a nontrivial Lie symmetry of a given  RDC equation  of the form (\ref{1-1-0}) are derived  by using DEs  and the group of ETs.  In Section 4, sufficient conditions are obtained and, as a result, the LSC problem is solved up to  the equivalence transformations. In Section 5,  form-preserving transformations (FPTs) are constructed. Applying these transformations to the RDC  equations derived in Section 4, we have shown that there are exactly 9 correctly-specified FPTs (they do not belong to ETs !), which allow us to map a RDC  equation with  nontrivial symmetry to another equation with the same symmetry.
 Finally, the main result of the paper is derived, namely: the  canonical list of all inequivalent RDC equations  admitting nontrivial Lie symmetry (up to any point  transformations)  is obtained.  In Section 6, an example devoted to the  construction of  exact solutions for the porous-Fisher type equation with the Burgers  term is presented.
 Finally, some discussion and conclusions are given in Section 7.

\section{\bf Determining equations for finding Lie symmetries and   the group of   equivalence transformations}

In general, the this  section  and Section 3, 4 and 5 are devoted to realization of the following algorithm for solving LSC problem for the class of RDC equations

\begin{equation}\label{l1}
u_{t}=(D(u)u_{x})_{x}+(D(u)u_{y})_{y}+K^{1}(u)u_{x}+K^{2}(u)u_{y}+R(u).
\end{equation}
The algorithm consists of the following steps

\begin{enumerate}
\item Application of  the classical Lie method for deriving the system of  DEs.
    \item  Finding the  principal algebra $A^{pr}$.
\item Construction of the group of equivalence transformations  (ETs)
    for the class  of RDC  equations in question.
 \item Deriving necessary conditions for possible extensions of $A^{pr}$ i.e. existence of nontrivial Lie symmetry.
     \item  Finding sufficient  conditions for  extensions of $A^{pr}$ and
deriving LSC for the  class of equations (\ref{l1}) using the group of ETs.
\item Construction of form-preserving transformations (FPTs).
\item Deriving a  complete  LSC  using   FPTs.
\end{enumerate}

Of course, this scheme can be modified  in some cases,
however, one may claim that it is a typical way to derive the so called canonical list of all inequivalent  equations from the class in question admitting nontrivial Lie symmetry (see more details in Chapter 2 of
\cite{ch-se-pl-2018}).

In this section,  we implement  the first, second and third steps, hence the following statement can be formulated.

\begin{theo}
The  principal algebra $A^{pr}$  of the class of RDC equations (\ref{l1})
is the three-dimensional Abelian algebra with the basic operators
\begin{equation}\label{m1}
A^{pr}=<\partial_{t}=\frac{\partial}{\partial t},\quad
\partial_{x}=\frac{\partial}{\partial x},\quad
\partial_{y}=\frac{\partial}{\partial y}>.
\end{equation}
\end{theo}

\textbf{Proof.}
To prove this statement we should construct  the system of DEs  for the class of equations in question. According to the classical Lie scheme (see e.g., \cite{ch-se-pl-2018, ol-86, fu-sh-se-93}), we write the most general structure of Lie symmetry operator (another terminology is infinitesimal operator) in the form
\begin{equation}\begin{array}{l}\label{m2}
X=\xi^{0}(t,x,y,u)\partial_{t}+\xi^{1}(t,x,y,u)\partial_{x}+\xi^{2}(t,x,y,u)\partial_{y}+\eta(t,x,y,u)\partial_{u},
\end{array}\end{equation}
where  $
\xi^{0},\;\xi^{1},\;\xi^{2},\;\eta$ are to-be-determined functions.
On the other hand, the class of equations
(\ref{l1}) is considered as the manifold
\be \label{2-4-4*}
{S} \equiv \Bigl\{(D(u)u_{x})_{x}+(D(u)u_{y})_{y}+K^{1}(u)u_{x}+K^{2}(u)u_{y}+R(u)- u_t =0 \Bigr\}\ee
in the prolonged space of the  variables
\[
t,\,  x,\, y, \, u,\,  \underset{1}{u},\,  \underset{2}{u}
\]
where $\underset{1}{u}\equiv(u_{t}, u_{x},u_{y})$ and  $\underset{2}{u}\equiv(u_{tt}, u_{tx},u_{ty},u_{xx},u_{xy},u_{yy})$.

An equation of the form (\ref{l1}) is invariant under the transformations generated by the
infinitesimal operator \eqref{m2}
when the following invariant criteria is satisfied:
\be \label{m2*}
 \underset{2}X \Bigl( (D(u)u_{x})_{x}+(D(u)u_{y})_{y}+K^{1}(u)u_{x}+K^{2}(u)u_{y}+R(u)- u_t =0\Bigr)\Bigl|_ {S}= 0,
\ee
where $ \underset{2}X$ is the second prolongation of $X$, which is calculated by the well-known formulae (see, e.g., \cite{ch-se-pl-2018, ol-86, fu-sh-se-93}). Making straightforward calculations (nowadays it can be done using computer algebra packages, e.g., Maple), one arrives at a system of DEs for finding the coefficients of the  infinitesimal operator \eqref{m2}. The system of DEs consists of the following differential equations:

\begin{equation}\begin{array}{l}\label{m3}
\xi^{0}_{x}=\xi^{0}_{y}=\xi^{0}_{u}=\xi^{1}_{u}=\xi^{2}_{u}=\eta_{uu}=0,
\end{array}\end{equation}
\begin{equation}\begin{array}{l}\label{m4}
\xi^{1}_{x}=\xi^{2}_{y},\quad \xi^{2}_{x}+\xi^{1}_{y}=0,
\end{array}\end{equation}
\begin{equation}\begin{array}{l}\label{m5}
\eta\dot{D}=(2\xi^{1}_{x}-\xi^{0}_{t})D,
\end{array}\end{equation}
\begin{equation}\begin{array}{l}\label{m6}
\eta\dot{K}^{1}=(\xi^{1}_{x}-\xi^{0}_{t})K^{1}-\xi^{2}_{x}
K^{2}-2\eta_{xu}D-
2\eta_{x}\dot{D}-\xi^{1}_{t},
\end{array}\end{equation}
\begin{equation}\begin{array}{l}\label{m612}
\eta\dot{K}^{2}=\xi^{2}_{x}
K^{1}+(\xi^{1}_{x}-\xi^{0}_{t})K^{2}-2\eta_{yu}D
-2\eta_{y}\dot{D}-\xi^{2}_{t},
\end{array}\end{equation}
\begin{equation}\begin{array}{l}\label{m7}
\eta\dot{R}=(a-\xi^{0}_{t})R-\triangle\eta D-\eta_{x}K^{1}-\eta_{y}K^{2}+\eta_{t},
\end{array}\end{equation}
 where the Laplacian $\Delta=\frac{\partial^{2}}{\partial x^{2}}+\frac{\partial^{2}}{\partial y^{2}}$  and the upper dot means differentiation w.r.t. $u$.
 The linear  equations  (\ref{m3})-(\ref{m4}) can be easily integrated. However, the general solution of (\ref{m5})-(\ref{m7}) depends essentially on the form of the functions $D(u), K^{1}(u), K^{2}(u)$ and $R(u)$. In the case of finding the  principal algebra, the problem simplifies because the functions $\xi^{0}$,\;$\xi^{1}$,\;$\xi^{2}$ and  $\eta$ should satisfy the system of DEs for arbitrary smooth functions
 $D(u), K^{1}(u), K^{2}(u)$ and $R(u)$. So, we immediately obtain $\eta=0$ from Eq. (\ref{m5}), while the system of linear equations
 \begin{equation}\label{m8}
\xi^{0}_{t}=\xi^{1}_{t}=\xi^{2}_{t}=\xi^{0}_{x}=\xi^{1}_{x}=
\xi^{2}_{x}=\xi^{0}_{y}=\xi^{1}_{y}=\xi^{2}_{y}=\xi^{0}_{u}=
\xi^{1}_{u}=\xi^{2}_{u}=0
\end{equation}
 is derived for the functions $\xi^{0}$,\;$\xi^{1}$  and $\xi^{2}$.
 As a result, the general solution of the  system of DEs  takes the form

\begin{equation}\begin{array}{l}\label{m9}
\xi^{0}=c_{0},\quad \xi^{1}=c_{1},\quad \xi^{2}=c_{2},\quad \eta=0,
\end{array}\end{equation}
 where $c_{0}, c_{1}$ and $c_{2}$  are arbitrary constants. The Lie symmetry operator  (\ref{m2}) with the coefficients
(\ref{m9}) produces the the three-dimensional Abelian algebra  (\ref{m1}).

The proof is now completed. \quad $\Box$

Now we turn to the group of  ETs. In order to find this group for the class of RDC equations (\ref{l1}),  the standard technique, which was formalized in \cite{akh91} (see also the relevant chapters in \cite{ch-se-pl-2018} and \cite{blu-ch-an}), can be  applied.

\begin{theo}
The group of the continuous ETs  of the  class  of RDC equations (\ref{l1})
is the 10-parameter Lie group

\begin{equation}\begin{array}{l}\label{k1}
t'=e^{\theta_{0}}t+m_{0},\\
x'=e^{\theta_{1}}(xcos\theta_{2}-ysin\theta_{2})+g_{1}t+m_{1},\\
y'=e^{\theta_{1}}(ycos\theta_{2}+xsin\theta_{2})+g_{2}t+m_{2},\\
u'=e^{\theta}u+m,\\
D'=e^{2\theta_{1}-\theta_{0}}D,\\
K^{1'}=e^{-\theta_{0}}[e^{\theta_{1}}(K^{1}cos\theta_{2}-K^{2}sin\theta_{2})-q_{1}],\\
K^{2'}=e^{-\theta_{0}}[e^{\theta_{1}}(K^{2}cos\theta_{2}+K^{1}sin\theta_{2})-q_{2}],\\
R'=e^{\theta-\theta_{0}}R,
\end{array}\end{equation}
 where \; $g_{1},\;g_{2},\;\theta_{0},\; \theta_{1},\; \theta_{2},\; \theta,\;  m_{0},\;
m_{1},\; m_{2},\;$  and $m$  are arbitrary group
parameters.

\end{theo}

\textbf{Proof} In order to find  operator $E$, we use the standard procedure (see, e.g., section 2.3.1 \cite{ch-se-pl-2018}), which is based on a
modification of the classical Lie method. In the case of class (\ref{l1}), we should start from the  infinitesimal operator

\begin{equation}\begin{array}{l}\label{k2}
E=\xi^{0}\partial_{t}+\xi^{1}\partial_{x}+\xi^{2}\partial_{y}+\eta\partial_{u}+\zeta^{0}\partial_{D}+\zeta^{1}\partial_{K^{1}}+
\zeta^{2}\partial_{K^{2}}+\zeta^{3}\partial_{R},
\end{array}\end{equation}
 where the coefficients
 \begin{equation}\begin{array}{l}\label{k2*}
 \medskip
 \xi^{0}=\xi^{t}(t,x,y,u),\quad \xi^{1}=\xi^{1}(t,x,y,u),\quad \xi^{2}=\xi^{2}(t,x,y,u),\quad \eta=\eta(t,x,y,u), \\
  \medskip
\zeta^{0}=\zeta^{0}(t,x,y,u,D,K^{1},K^{2},R), \quad
\zeta^{1}=\zeta^{1}(t,x,y,u,D,K^{1},K^{2},R),
\\
\zeta^{2}=\zeta^{2}(t,x,y,u,D,K^{1},K^{2},R), \quad \zeta^{3}=\zeta^{3}(t,x,y,u,D,K^{1},K^{2},R)
\end{array}\end{equation}
are to-be-determined functions. Note that the coefficients $\zeta^{j}, j=0,...,3$ may depend on $D, K^{1}, K^{2}$ and/or $R$ (in contrast to other coefficients of $E$).

  Lie's invariance criterion should be applied to system of equations consisting of (\ref{l1}) and a set of differential consequences of the functions $D$, $K^{1}$, $K^{2}$ and  $R$ with respect to the variables $t$, $x$, $y$, $u_{t}$, $u_{x}$ and $u_{y}$.
  Making relevant calculations, we arrive at the following system
\begin{equation}\begin{array}{l}\label{k4}
\xi^{0}_{u}=\xi^{1}_{u}=\xi^{2}_{u}=\xi^{0}_{x}=\xi^{0}_{y}=\eta_{uu}=0,
\end{array}\end{equation}
\begin{equation}\begin{array}{l}\label{k411}
\zeta^{0}_{t}=\zeta^{0}_{x}=\zeta^{0}_{y}=\zeta^{1}_{t}=\zeta^{1}_{x}=\zeta^{1}_{x}=\zeta^{2}_{t}=\zeta^{2}_{x}=\zeta^{2}_{x}=
\zeta^{3}_{t}=\zeta^{3}_{x}=\zeta^{3}_{y}=
\\
=\zeta^{0}_{u_{t}}=\zeta^{0}_{u_{x}}=\zeta^{0}_{u_{y}}=
\zeta^{1}_{u_{t}}=\zeta^{1}_{u_{x}}=\zeta^{1}_{u_{y}}=\zeta^{2}_{u_{t}}=\zeta^{2}_{u_{x}}=\zeta^{2}_{u_{y}}=\zeta^{2}_{u_{t}}=\zeta^{3}_{u_{x}}=\zeta^{3}_{u_{y}}
\end{array}\end{equation}
\begin{equation}\begin{array}{l}\label{k5}
\xi^{1}_{x}=\xi^{2}_{y},\quad  \xi^{2}_{x}+\xi^{1}_{y}=0,
\end{array}\end{equation}
\begin{equation}\begin{array}{l}\label{k6}
\zeta^{0}=(\xi^{0}_{t}-2\xi^{1}_{x})D,
\end{array}\end{equation}
\begin{equation}\begin{array}{l}\label{k71}
\zeta^{1}=\xi^{1}_{t}D-\xi^{1}_{x}K^{1}+\xi^{1}_{y}K^{2},
\end{array}\end{equation}
\begin{equation}\begin{array}{l}\label{k72}
\zeta^{2}=\xi^{2}_{t}D+\xi^{2}_{x}K^{1}+(\xi^{2}_{y}-2\xi^{1}_{x})K^{2},
\end{array}\end{equation}
\begin{equation}\begin{array}{l}\label{k8}
\zeta^{3}=(\eta_{u}-2\xi^{1}_{1})R
\end{array}\end{equation}
to find the coefficients $\xi^{0}$,\;  $\xi^{a}$,\;  $\eta$,\;  $\zeta^{0}$,\;  $\zeta^{1}$,\; $\zeta^{2}$,\;
$\zeta^{3}$ of operator (\ref{k2}).

Because the system of PDEs (\ref{k4})-(\ref{k8}) is linear, its general solutions can be easily constructed, hence we obtain
\begin{equation}\begin{array}{l}\label{k9}
\xi^{0}=\kappa_{0}t+d_{0},\\
\xi^{1}=\kappa_{1}x+g_{1}t-cy+d_{1},\\
\xi^{2}=\kappa_{1}y+g_{2}t+cx+d_{2},\\
\eta=\kappa u+d,\\
\zeta^{0}=(\kappa_{0}-2\kappa_{1})D,\\
\zeta^{1}=g_{1}D+(c-\kappa_{1})K^{1}-cK^{2},\\
\zeta^{2}=g_{2}D+cK^{1}+(c-\kappa_{1})K^{2},\\
\zeta^{3}=(\kappa-2\kappa_{1})R.
\end{array}\end{equation}
where $c,  \; \kappa, \; \kappa_{0}, \; \kappa_{1},  \; g_{1}, \; g_{2},\;
d_{0}, \; d_{1}, \; d_{2}$  and $d$ are arbitrary parameters.

Thus, the infinitesimal operator (\ref{k2}) with the coefficients from (\ref{k9}) generates the 10-dimensional Lie algebra. Finally, making  standard calculations one can easily show that this algebra   produces the 10-parameter Lie group of equivalence transformations (\ref{k1}).

The proof is now completed.\quad $\Box$

\section{\bf Necessary conditions for   existence of a nontrivial Lie symmetry}

Here we are searching for  necessary conditions, which are needed for  extension of the   principal Lie algebra
(\ref{m1}). In other words, we need to establish all possible forms  of the functions $D$, $K^{1}$, $K^{2}$,  and  $R$
 leading to extension of Lie symmetry of the relevant RDC equations from class (\ref{l1}).
 The main result of this section can be formulated as follows.

\begin{theo}\label{T2-5-2} If  an arbitrary equation belonging to the class of RDC equations (\ref{l1})
admits a  maximal algebra of invariance (MAI) of a higher  dimensionality than algebra   (\ref{m1}), then the functions $D$, $K^{1}$, $K^{2}$,and  $R$ must possess the structures listed in Table  \ref{Tabl1-1}, where  $\lambda_{3},\; \lambda_{4},\; \lambda_{5},\; k,\; m,\; p$  and $s$ are arbitrary constants. Any other RDC equation possessing a nontrivial Lie symmetry is reducible to one of those from Table  \ref{Tabl1-1} by an appropriate ET of the form (\ref{k1}).
\end{theo}

{\begin{longtable}[c]{|c|c|c|c|c|c|}
\caption[c]{Necessary conditions for   nontrivial Lie symmetries  }\label{Tabl1-1}\\
\hline
 & $D$ & $K^{1}$ & $K^{2}$ & $R$ & Restrictions \\
\hline 1 & $\forall$ & $0$& $0$&  $\forall$&\\
\hline 2 & $e^{s u}$ & $e^{mu}\cos(pu)$ & $e^{mu}\sin(pu)$ &
$\lambda_{3}e^{(2m-s)u}$ &$(m,p) \neq (0,0), m\neq s$\\
\hline 3& $e^{u}$&$u$& $0$& $\lambda_{3}e^{-u}$ &\\
\hline 4 & $e^{u}$ & $e^{u}\cos(pu)$ & $e^{u}\sin(pu)$ &
$\lambda_{3}e^{u}+\lambda_{4}$ &\\
\hline 5& $u^{k}$&$u^{m}\cos(p\ln u)$& $u^{m}\sin(p\ln u)$&
$\lambda_{3}u^{2m-k+1}$ &$m\neq k, (m,p)\neq (0,0)$\\
\hline 6& $u^{k}$&$\ln u$& $0$&
$\lambda_{3}u^{-k+1}$ &$k\neq 0$\\
\hline 7& $u^{k}$&$u^{k}\cos(p\ln u)$& $u^{k}\sin(p\ln u)$&
$u(\lambda_{3}u^{k}+\lambda_{4})$ &$k\neq 0$\\
\hline 8 & $1$ & $u$& $0$& $\lambda_{3}u+\lambda_{4}$&\\
\hline 9 & $1$ & $\ln u$& $0$& $u(\lambda_{3} \ln^{2}
u+\lambda_{4}\ln u+\lambda_{5})$&\\
\hline
\end{longtable}

\textbf{Proof} To prove this theorem, one needs to analyze the system of DEs (\ref{m3})-(\ref{m7}). Obviously, Eqs. (\ref{m3})-(\ref{m4}) are rather simple and can be easily integrated. In particular, Eqs. (\ref{m3}) allow us to obtain
\be\label{3-1}  \xi^{0}=\xi^{0}(t),  \quad \xi^{i}=\xi^{i}(t,x,y), \quad \eta= a(t,x,y)u+b(t,x,y),
\ee
where $a(t,x,y)$  and $b(t,x,y)$ are arbitrary functions, $i=1,2$.

 The system of equations  (\ref{m5})-(\ref{m7})
 from the formal point of view it is more complicated object than the general RDC equation  (\ref{l1}). However, unknown functions in Eqs. (\ref{m5})-(\ref{m7})  depend on different variables (e.g., $\xi^0$ depends on $t$ while $D$ depends on $u$ only)  and it is
  a common peculiarity  of such type  systems, which allows us to work out an algorithm for their solving.
  Unfortunately, this algorithm usually is quite cumbersome and consists of examination of several inequivalent cases. Happily we can partly use the algorithm presented in Chapter 2 \cite{ch-se-pl-2018} for the solving similar system obtained for (1+1)-dimensional general RDC equation.
In particular, to simplify the relevant calculations we introduce the so called structural constants as follows

\begin{equation}\begin{array}{l}\label{p1}
a=k_{1}\varphi,\;  b=k_{2}\varphi,  \; 2\xi^{1}_{x}-\xi^{0}_{t}=k\varphi,   \; \xi^{1}_{x}-\xi^{0}_{t}=m\varphi,  \; \xi^{2}_{x}=p\varphi, \; a_{x}=\alpha_{1}\varphi,\; a_{y}=\alpha_{2}\varphi,\; \\
b_{x}=\beta_{1}\varphi,\;  b_{y}=\beta_{2}\varphi,\; \xi^{1}_{t}=\gamma_{1}\varphi,\; \xi^{2}_{t}=\gamma_{2}\varphi,\;  a-\xi^{0}_{t}=(2m+k_{1}-k)\varphi,\; \triangle a=h_{1}\varphi, \;
\\ \triangle b=n_{1}\varphi,\; a_{t}=h_{2}\varphi, \;  b_{t}=n_{2}\varphi,\;
\end{array}\end{equation}
where
$k,\;m,\;p,\;k_{i},\;n_{i},\;h_{i},\;\alpha_{i},\;\beta_{i}$  and $\gamma_{i} \, (i= 1,2)$
are some structural constants relating all unknown functions with the function
$\varphi(t,x,y)$, which is arbitrary at the moment.
Using notations  (\ref{p1}),  Eqs. (\ref{m5})-(\ref{m7}) can be rewritten in the form
\begin{equation}\begin{array}{l}\label{p2*}
(k_{1}u+k_{2})\dot{D}=kD, \\
(k_{1}u+k_{2})\dot{K^{1}}=mK^{1}-pK^{2}-2\alpha_{1}D-2(\alpha_{1}u+\beta_{1})\dot{D}-\gamma_{1},\\
(k_{1}u+k_{2})\dot{K^{2}}=pK^{1}+
mK^{2}-2\alpha_{2}D-2(\alpha_{2}u+\beta_{2})\dot{D}-\gamma_{2},\\
(k_{1}u+k_{2})\dot{R}=(2m-k+k_{1})R-(h_{1}u+n_{1})D-
(\alpha_{1}u+\beta_{1})K^{1}-
\\
-(\alpha_{2}u+\beta_{2})K^{2}+h_{2}u+n_{2}.
\end{array}\end{equation}
 System (\ref{p2*})  possesses a simpler structure comparing with Eqs. (\ref{m5})-(\ref{m7}) because one does not involve the functions on the variables $t,\,  x$  and $ y$. On the other hand, this system produces all possible forms the functions $D$, $K^{1}$, $K^{2}$  and  $R$  leading to extensions  of the   principal Lie algebra (\ref{m1}). The rest of the proof is devoted to solving system (\ref{p2*}), which consists of four linear ODEs.

First of  all, it can be noted that system (\ref{p2*}) can be slightly simplified using the correctly-specified ETs of the form
\begin{equation}\begin{array}{l}\label{p3}
t'=t,\; x'=x+\theta_{1}t,\; y'=y+\theta_{2}t,\; u'=u,\;
D'=D,\; K^{1'}=K^{1}-\theta_{1},\;
\\
K^{2'}=K^{2}-\theta_{2},\; R'=R,
\end{array}\end{equation}
which are a particular case of (\ref{k1}).
In fact, the parameters $\theta_{1}$  and $ \theta_{2}$  can be chosen in a such way that system (\ref{p2*})  can be transformed to the same form with $\gamma_{1}=\gamma_{2}=0$ (hereafter primes are skipped)
\begin{equation}\begin{array}{l}\label{p3}
(k_{1}u+k_{2})\dot{D}=kD, \\
(k_{1}u+k_{2})\dot{K^{1}}=mK^{1}-pK^{2}-2\alpha_{1}D-2(\alpha_{1}u+\beta_{1})\dot{D},\\
(k_{1}u+k_{2})\dot{K^{2}}=pK^{1}+
mK^{2}-2\alpha_{2}D-2(\alpha_{2}u+\beta_{2})\dot{D},\\
(k_{1}u+k_{2})\dot{R}=(2m-k+k_{1})R-(h_{1}u+n_{1})D-
(\alpha_{1}u+\beta_{1})K^{1}-
\\
-(\alpha_{2}u+\beta_{2})K^{2}+h_{2}u+n_{2},
\end{array}\end{equation}
provided
\begin{equation}\label{p7}
m^{2}+p^{2}\neq 0,
\end{equation}
The possibility $m=p=0$ will be treated below when one comes up.

It can be noted that the first equation in (\ref{p3}) has the same structure as one in the system of DEs for the class of RDC equations (\ref{1-1-0}) with $n=1$ (see Eq.(2.80) in \cite{ch-se-pl-2018}). Thus, the following five different case should be examined  (see P.41 in \cite{ch-se-pl-2018}):

1)$k_{1}=0,\; k_{2}=0,\; k=0$;

2)$k_{1}=0,\; k_{2}\neq0,\;k\neq0$;

3)$k_{1}\neq0,\; k_{2}=0,\; k\neq0$;

4)$k_{1}=0,\; k_{2}\neq 0,\; k=0$;

5)$k_{1}\neq 0,\; k_{2}=0,\; k=0$.

Since Case 1)  is rather trivial (the first equation in (\ref{p3}) simply vanishes) and leads to the first case of Table \ref{Tabl1-1}. Here the relevant analysis is omitted.

Consider Case 2), i.e. $k_{1}=0,\; k_{2}\neq 0,\; k\neq0$ (without loosing a generality we can set
$k_{2}=k=1$). The first two equations of system  (\ref{p1}) immediately give
$a=0,\;b=\varphi$, hence
\begin{equation}\begin{array}{l}\label{p8}
\alpha_{1}=\alpha_{2}=h_{1}=h_{2}=0
\end{array}\end{equation}
and
\begin{equation}\begin{array}{l}\label{p9}
b=2\xi^{1}_{x}-\xi^{0}_{t},   \; mb=\xi^{1}_{x}-\xi^{0}_{t},  \; pb=\xi^{2}_{x}, \; \gamma_{1}b=\xi^{1}_{t}, \; \gamma_{2}b=\xi^{2}_{t}.
\end{array}\end{equation}

It follows from  Eqs. (\ref{p9}) that
\begin{equation}\begin{array}{l}\label{p10}
(2m-1)\xi^{1}_{x}=(m-1)\xi^{0}_{t}.
\end{array}\end{equation}
Now we take differential consequences of (\ref{p10}) w.r.t.$x$ and  $y$
and use Eq. (\ref{m4}). As a result, the equations
\begin{equation}\begin{array}{l}\label{p12}
m(2m-1)b_{x}=0,  \; m(2m-1)b_{y}=0.
\end{array}\end{equation}
are derived.
Finally, differentiating the last equation in (\ref{p9}) w.r.t.$x$ and  $y$,  and using Eqs. (\ref{p12}), the compatibility constrain
\begin{equation}\begin{array}{l}\label{p14}
m(m-1)(2m-1)b_{t}=0
\end{array}\end{equation}
is obtained.

The compatibility  constraint (\ref{p12}) leads to four different subcases, which  can be considered
step by step.
Namely, the following  different subcases should be examined:  (2i) $m\neq0,\ \frac{1}{2},\ k$, (2ii) $m=0$,
 (2iii) $m=\frac{1}{2}$  and (2iv) $m=1$.

We start from the most general one (2i).
Obviously  Eqs. (\ref{p12}) and  (\ref{p14}) immediately give   $b=const$, hence
\begin{equation}\begin{array}{l}\label{p15}
\beta_{1}=\beta_{2}=n_{1}=n_{2}=0.
\end{array}\end{equation}

Thus, system  (\ref{p2*}) takes the form
\begin{equation}\begin{array}{l}\label{p16}
\dot{D}=D,\quad
\dot{K^{1}}=mK^{1}-pK^{2},\quad
\dot{K^{2}}=pK^{1}+mK^{2},\quad
\dot{R}=(2m-1)R.
\end{array}\end{equation}
The general solution of the latter is
\begin{equation}\begin{array}{l}\label{p160}
D=\lambda_{0}e^{u},\quad
K^{1}=\lambda_{1}e^{mu}\cos(pu+\lambda_{2}),\quad
K^{2}=\lambda_{1}e^{mu}\sin(pu+\lambda_{2}),\quad
R=\lambda_{3}e^{(2m-1)u}.
\end{array}\end{equation}
Hereafter  $\lambda_{0}\neq 0,\; \lambda_{1}\neq 0,\;\lambda_{2},\; \lambda_{3}$
are arbitrary constants.
However, three of these constants can be reduced to $\lambda_{0}=\lambda_{1}=1,\;\lambda_{2}=0$ using the equivalence transformation
\begin{equation}\begin{array}{l}\label{p180}
t\rightarrow\theta_{0}t,\quad
x\rightarrow\theta_{1}x,\quad
y\rightarrow\theta_{1}y,\quad
u\rightarrow u+\theta_{2}.
\end{array}\end{equation}
where
\begin{gather*}
\theta_{0}=\frac{\lambda_{0}e^{(2m-1)\frac{\lambda_{2}}{p}}}{\lambda^{2}_{1}}, \quad \theta_{0}=\frac{\lambda_{0}e^{(m-1)\frac{\lambda_{2}}{p}}}{\lambda_{1}},\quad \theta_{0}=-\frac{\lambda_{2}}{p}.
\end{gather*}

Thus, the second case with $s=1$ of   Table \ref{Tabl1-1} is identified.

Consider subcase  (2ii) $m=0$. In a quite similar way as it was done in case (2i), using  system (\ref{p1}), one may extract   the restrictions
\begin{equation}\begin{array}{l}\label{p17}
\beta_{1}=\beta_{2}=n_{1}=n_{2}=0.
\end{array}\end{equation}
Thus, system  (\ref{p2*}) takes the form
\begin{equation}\begin{array}{l}\label{p18}
\dot{D}=D,\quad
\dot{K^{1}}=-pK^{2}-\gamma_{1},\quad
\dot{K^{2}}=pK^{1}-\gamma_{2},\quad
\dot{R}=-R.
\end{array}\end{equation}

The general solution of (\ref{p18}) depends on the constant $p$.
Assuming  $p\neq 0$,
we arrive at the following  general solution of (\ref{p18})
\begin{gather*}
D=\lambda_{0}e^{u},\quad
K^{1}=\lambda_{1}\cos(pu+\lambda_{2}),\quad
K^{2}=\lambda_{1}\sin(pu+\lambda_{2}),\quad
R=\lambda_{3}e^{-u},
\end{gather*}
So, the formulae (\ref{p160})
are valid also for  $m=0$.
Moreover, using ET (\ref{p180}) the coefficients in the above formulae are reducible  to
 $\lambda_{0}=\lambda_{1}=1,\;\lambda_{2}=0$.

Assuming  $p= 0$, we note that  restriction (\ref{p7}) is broken, hence  the constants $\gamma_{1}$  and $\gamma_{2}$ can be non-zero in (\ref{p18}). So, the  general solution of (\ref{p18}) is
\begin{gather*}
D=\lambda_{0}e^{u}, \quad K^{1}=\lambda_{1}u,\quad K^{2}=\lambda_{2}u,\quad
R=\lambda_{3}e^{-u}.
\end{gather*}
Now we again use the following ET
\begin{equation}\begin{array}{l}\label{p184}
t\rightarrow\frac{1}{\lambda_{0}}t,\quad
x\rightarrow\frac{\lambda_{0}}{\vec{\lambda}^{2}}(\lambda_{1}x+\lambda_{2}y),\quad
y\rightarrow\frac{\lambda_{0}}{\vec{\lambda}^{2}}(\lambda_{1}y-\lambda_{2}x),\quad
u\rightarrow u.
\end{array}\end{equation}
in order to simplify the coefficients as follows $\lambda_{0}=\lambda_{1}=1,\;\lambda_{2}=0$.

Thus, the third case  of   Table \ref{Tabl1-1} is identified.

Consider subcase  (2iii) $m=\frac{1}{2}$.
In a quite similar way  the corresponding system and its    general solution

\[D=e^{u},\quad
K^{1}=e^{\frac{1}{2}u}\cos pu,\quad
K^{2}=e^{\frac{1}{2}u}\sin pu,\quad
R=\lambda_{3} \]
were constructed.
So, formulae  (\ref{p160}) are valid also for  $m=\frac{1}{2}$.

Finally, subcase  (2iv) $m=1$ was examined.
It was shown that  system  (\ref{p2*}) takes the form
\begin{equation}\begin{array}{l}\label{p20}
\dot{D}=D,\quad
\dot{K^{1}}=K^{1}-pK^{2},\quad
\dot{K^{2}}=K^{2}+pK^{1},\quad
\dot{R}=R+q_{2}.
\end{array}\end{equation}
The latter   possesses the general solution
\begin{gather*}
D=\lambda_{0}e^{u},\quad
K^{1}=\lambda_{1}e^{u}\cos(pu+\lambda_{2}),\quad
K^{2}=\lambda_{1}e^{u}\sin(pu+\lambda_{2}),\quad
R=\lambda_{3}e^{u}+\lambda_{4},
\end{gather*}
 where $\lambda_{0},\; \lambda_{1},\; \lambda_{2},\;\lambda_{3},\;\lambda_{4}$ are arbitrary constants.
Here  three coefficients  are again reducible  to
 $\lambda_{0}=\lambda_{1}=1,\;\lambda_{2}=0$  via application of  ET (\ref{p180}).

As a result, the fourth case  of   Table \ref{Tabl1-1} is identified.

Thus, the cases 1, 2 (with $s\not=0$), 3 and  4   of  Table \ref{Tabl1-1} are identified.
 All other cases of  Table \ref{Tabl1-1} were  derived in a very similar way by the examination of Cases 3)---5).

 The proof is now completed.\quad $\Box$

 \newpage

 \section {\bf Lie symmetry classification   using the equivalence transformations}

 As one may note, the steps 1--4 of the algorithm presented in Section 2 are already realized. Here we are going to identify  sufficient conditions needed for  extension of the   principal Lie algebra (\ref{m1}) (see step 5 of the algorithm). Having this done and taking into account  the group of ETs (\ref{k1}), we can complete LSC for the  class of  RDC equations (\ref{l1}).

\begin{theo}\label{T2-5-4}
All possible nontrivial  MAI (i.e., Lie algebras  of dimensionality four and higher)
of  RDC equations  of the form \eqref{l1}  depending on the
functions   $D,\,  K^1, \, K^1$  and $R$
are presented in  Table \ref{Tabl2-2}.  Any  other equation of the form \eqref{l1} with a nontrivial Lie symmetry  is reduced by an ET from $\mathcal{E}$ (\ref{k1}) to one of 32 equations  listed  in Table \ref{Tabl2-2}.

\end{theo}

\begin{longtable}[c]{|c|c|c|c|}
\caption [c]{LSC of the  class of  RDC equations (\ref{l1}) using the group of  ETs (\ref{k1}) }\label{Tabl2-2}\\
\hline
 & Equation  & MAI &  Restrictions \\
\hline 1 & $u_{t}=(D(u)u_{x})_{x}+(D(u)u_{y})_{y}+$& $<\partial_{t},
\partial_{x}, \partial_{y}, J_{12}>$ & $D-\forall, R-\forall$ \\
&$+R(u)$& $$ & $$ \\
\hline 2 & $u_{t}=(D(u)u_{x})_{x}+(D(u)u_{y})_{y}$ & $<\partial_{t},
\partial_{x}, \partial_{y}, J_{12}, D_{0}>$ & $D-\forall$ \\
\hline 3 & $u_{t}=\triangle u$ & $<\partial_{t},
\partial_{x}, \partial_{y}, J_{12}, G_{x}, G_{y}, I,$ & $$ \\
&& $D_{0}, \Pi,
 Q^{1}_{\infty}>$ & $$
\\
\hline 4 & $u_{t}=\triangle u +\gamma_{1}$ & $<\partial_{t},
\partial_{x}, \partial_{y}, J_{12}, G_{x}+\frac{1}{2}\gamma_{1}tx\partial_{u},$ & $$ \\ [0.5mm]
&& $G_{y}+\frac{1}{2}\gamma_{1}ty\partial_{u}, I-\gamma_{1}t\partial_{u},D_{0}+2\gamma_{1}t\partial_{u}, $ & $$ \\ [0.5mm]
&& $\Pi+\gamma_{1}t(2t+\frac{x^{2}+y^{2}}{4})\partial_{u},  Q^{1}_{\infty}>$ & $$\\ [0.5mm]
\hline 5 & $u_{t}=\triangle u +\gamma_{1}u$ & $<\partial_{t},
\partial_{x}, \partial_{y}, J_{12}, G_{x}, G_{y}, I, $ & $$ \\ [0.5mm]
&& $D_{0}+2\gamma_{1}tu\partial_{u},\Pi+\gamma_{1}t^{2}u\partial_{u}, Q^{2}_{\infty}> $ & $$ \\ [0.5mm]
\hline 6 & $u_{t}=\triangle u+\gamma_{1}u \ln u$ & $<\partial_{t},
\partial_{x}, \partial_{y},  J_{12}, e^{\gamma_{1} t}I, {\cal G}_{x}, {\cal G}_{y}>$ & $$
\\ [0.5mm]
\hline 7 & $u_{t}=(e^{u}u_{x})_{x}+(e^{u}u_{y})_{y}$ & $<\partial_{t},
\partial_{x}, \partial_{y}, J_{12}, D_{0}, D_{2}>$ & $\delta=1$
\\ [0.5mm]
\hline 8 & $u_{t}=(e^{\delta u}u_{x})_{x}+(e^{\delta u}u_{y})_{y}+\gamma_{1}e^{mu}$ & $<\partial_{t},
\partial_{x}, \partial_{y}, J_{12}, $ & $m \neq 0$ \\
&& $(\delta-m)D_{0}-2 D_{4}>$ & $$
\\
\hline 9 & $u_{t}=(e^{u}u_{x})_{x}+(e^{u}u_{y})_{y} +\gamma_{1}$ & $<\partial_{t},
\partial_{x}, \partial_{y}, J_{12}, D_{2},T_{2}>$ & $\delta=1$
\\
\hline 10 & $u_{t}=(e^{u}u_{x})_{x}+(e^{u}u_{y})_{y}+\gamma_{1}e^{u}+\gamma_{2}$ & $<\partial_{t},
\partial_{x}, \partial_{y}, J_{12}, T_{2}>$ & $$ \\
\hline 11 & $u_{t}=(u^{k}u_{x})_{x}+(u^{k}u_{y})_{y}$ & $<\partial_{t},
\partial_{x}, \partial_{y}, J_{12}, D_{0}, D_{1}>$ & $k \neq -1; 0$\\ [0.5mm]
\hline 12 & $u_{t}=(u^{k}u_{x})_{x}+(u^{k}u_{y})_{y}+\gamma_{1}u^{m}$ & $<\partial_{t},
\partial_{x}, \partial_{y}, J_{12}, $ & $k \neq -1; 0$ \\
&& $(m-1)D_{0}-D_{1}>$
&$m \neq 1$\\ [0.5mm]
\hline 13 & $u_{t}=(u^{k}u_{x})_{x}+(u^{k}u_{y})_{y}+\gamma_{1}u$ & $<\partial_{t},
\partial_{x}, \partial_{y}, J_{12}, D_{1},T_{1}>$ & $k \neq -1; 0$
\\ [0.5mm]
\hline 14 & $u_{t}=(u^{k}u_{x})_{x}+(u^{k}u_{y})_{y}+$ & $<\partial_{t},
\partial_{x}, \partial_{y}, J_{12}, T_{1}>$ & $k \neq 0$
\\ [0.5mm]
& $+\gamma_{1}u^{k+1}+\gamma_{2}u$ & $$& $$
\\ [0.5mm]
\hline 15 & $u_{t}=(u^{-1}u_{x})_{x}+(u^{-1}u_{y})_{y}$ & $<\partial_{t},
\partial_{x}, \partial_{y}, J_{12}, D_{3},X_{\infty}>$ & $k=-1$
\\ [0.5mm]
\hline 16 & $u_{t}=(u^{-1}u_{x})_{x}+(u^{-1}u_{y})_{y}+\gamma_{1}u$ & $<\partial_{t},
\partial_{x}, \partial_{y}, J_{12}, T_{1}, X_{\infty}>$ & $k=-1$ \\ [1mm]
\hline 17& $u_{t}=(e^{\delta
u}u_{x})_{x}+(e^{\delta
u}u_{y})_{y}+$&
$<\partial_{t}, \partial_{x},\partial_{y},$&$m\neq \delta$\\ [1mm]
&$+e^{mu}[u_{x}\cos(pu)+u_{y}\sin(pu)]+$
&$(m-\delta)D_{0}+D_{4}+pJ_{12}>$&$(m,p)\neq (0,0)$\\ [1mm]
&$+\sigma e^{(2m-\delta)u}$
&$$&$$\\ [0.5mm]
\hline 18 & $u_{t}=(e^{u}u_{x})_{x}+(e^{u}u_{y})_{y}+$ & $<\partial_{t}, \partial_{x},\partial_{y}, $&$\delta=1$\\ [0.5mm]
&$+uu_{x}+\sigma e^{-u}$&$D_{0}-D_{4}-t\partial_{x}>$&$$
\\ [0.5mm]
\hline 19 & $u_{t}=(e^{u}u_{x})_{x}+(e^{u}u_{y})_{y}+$ & $<\partial_{t}, \partial_{x}, \partial_{y},
D_{4}+pJ_{12}>$ & $\delta=1$
\\ [0.5mm]
&$+e^{u}[u_{x}\cos(pu)+u_{y}\sin(pu)]+\sigma e^{u}$&$$& $$
\\ [0.5mm]
\hline 20 & $u_{t}=(e^{u}u_{x})_{x}+(e^{u}u_{y})_{y}+$ & $<\partial_{t}, \partial_{x},\partial_{y},
T_{2}>$& $$
\\ [0.5mm]
$$ & $+e^{u}u_{x}+\sigma e^{u}+\gamma_{1}$& $$ & \\ [0.5mm]
\hline 21&
$u_{t}=(u^{k}u_{x})_{x}+(u^{k}u_{y})_{y}+$&$<\partial_{t}, \partial_{x}, \partial_{y}, (m-k)D_{0}+$&
$m\neq k, $\\ [0.5mm]
&$+u^{m}[u_{x}\cos(p\ln u)+$&$+D_{3}+pJ_{12}>$&$(m,p)\neq (0,0)$
\\ [0.5mm]
&$+u_{y}\sin(p\ln u)]+\sigma u^{2m-k+1}$&$$&$$
\\ [0.5mm]
\hline 22& $u_{t}=(u^{k}u_{x})_{x}+(u^{k}u_{y})_{y}+$&$<\partial_{t}, \partial_{x}, \partial_{y}, $&$k\neq0$\\ [0.5mm]
&$+\ln uu_{x}+\sigma u^{-k+1}$&$kD_{0}-D_{3}-t\partial_{x}>$&$$
\\ [0.5mm]
\hline 23&
$u_{t}=(u^{k}u_{x})_{x}+(u^{k}u_{y})_{y}+$ & $<\partial_{t}, \partial_{x}, \partial_{y},D_{3}+pJ_{12}>
$& $k\neq0, p\neq0$\\ [0.5mm]
&$+u^{k}[u_{x}\cos(p\ln u)+$& &
\\ [0.5mm]
&$+u_{y}\sin(p\ln u)]+\sigma u^{k+1}$&$$&$$
\\ [0.5mm]
\hline 24& $u_{t}=(u^{k}u_{x})_{x}+(u^{k}u_{y})_{y}+$ & $<\partial_{t}, \partial_{x}, \partial_{y}, D_{3}>$&
$k\neq0$,\\ [0.5mm]
&$+u^{k}u_{x}+\lambda_{3} u^{k+1}$&&$k^{2}\lambda_{3}\neq 4(k+1)$
\\ [0.5mm]
\hline 25& $u_{t}=(u^{k}u_{x})_{x}+(u^{k}u_{y})_{y}+u^{k}u_{x}+$ & $<\partial_{t}, \partial_{x}, \partial_{y}, T_{1}> $&
$k\neq0$,\\ [0.5mm]
&$+\lambda_{3} u^{k+1}+\gamma_{1}u$&&$k^{2}\lambda_{3}\neq 4(k+1)$
\\ [0.5mm]
\hline 26& $u_{t}=(u^{k}u_{x})_{x}+(u^{k}u_{y})_{y}+$ & $<\partial_{t}, \partial_{x}, \partial_{y}, D_{3}, R_{1}, R_{2}>$& $k\neq -1; 0$
\\ [1mm]
&$+4\frac{k+1}{k}u^{k}u_{x}+4\frac{k+1}{k^{2}}u^{k+1}$&$$& $$
\\ [1mm]
\hline 27& $u_{t}=(u^{k}u_{x})_{x}+(u^{k}u_{y})_{y}+$&$<\partial_{t}, \partial_{x}, \partial_{y}, T_{1}, R_{1}, R_{2}>$& $k\neq -1; 0$\\ [1mm]
&$+4\frac{k+1}{k}u^{k}u_{x}+4\frac{k+1}{k^{2}}u^{k+1}+\gamma_{1}u$&$$&$$\\ [1mm]
\hline 28 & $u_{t}=\triangle u+uu_{x}+\gamma_{1}u$ & $<\partial_{t}, \partial_{x}, \partial_{y}, {\cal G}_{1}>$ & $$ \\ [0.5mm]
\hline 29 & $u_{t}=\triangle u+uu_{x}+\sigma$ & $<\partial_{t},
\partial_{x}, \partial_{y}, D_{0}-u\partial_{u}-$& $$\\ [0.5mm]
&& $-\frac{3}{2}\sigma t(G_{0}-\partial_{u}), G_{0}>$& $$
\\ [0.5mm]
\hline 30& $u_{t}=\triangle u+\ln uu_{x}+\sigma u$ & $<\partial_{t}, \partial_{x}, \partial_{y}, G_{1}>$&
$$  \\ [0.5mm]
\hline 31 & $u_{t}=\triangle u+\ln uu_{x}+\gamma_{1}u\ln u$ & $<\partial_{t}, \partial_{x}, \partial_{y},  {\cal G}_{2}>$&$$
\\ [0.5mm]
\hline 32 & $u_{t}=\triangle u+2\gamma_{1}\ln uu_{x}+$ &
$<\partial_{t},
\partial_{x}, \partial_{y},Y>$& \\ [0.5mm]
&$+u(\ln^{2} u+q)$&$$&$$
\\ [0.5mm]
\hline
\end{longtable}

\begin{remark} In Table \ref{Tabl2-2}, $k, \; m, \; p$  and $q,$ are arbitrary constants, $\sigma\in\{-1, 0, 1\}$, $\delta\in\{0,1\}$, $\gamma_{1},\;\gamma_{2}\in\{-1, 1\}$  and  the following designations for Lie symmetry operators  are introduced: \\ \medskip
$J_{12}=y\partial_{x}-x\partial_{y},\quad
D_{0}=2t\partial_{t}+x\partial_{x}+y\partial_{y},
\\
D_{1}=k x\partial_{x}+k y\partial_{y}+2u\partial_{u},\quad D_{2}=x\partial_{x}+ y\partial_{y}+2\partial_{u},
\\
D_{3}=k t\partial_{t}-u\partial_{u},\quad D_{4}=\delta t\partial_{t}-\partial_{u},
\\
G_{x}=t\partial_{x}-\frac{1}{2}xI,\quad G_{y}=t\partial_{y}-\frac{1}{2}yI,\quad I=u\partial_{u},
\\
\Pi=t^{2}\partial_{t}+tx\partial_{x}+ty\partial_{y}-(t+\frac{x^{2}+y^{2}}{4})u\partial_{u},
\\
{\cal G}_{x}=e^{\gamma_{1}t}(\partial_{x}-\frac{1}{2}\gamma_{1}xu\partial_{u}), \quad {\cal G}_{y}=e^{\gamma_{1}t}(\partial_{y}-\frac{1}{2}\gamma_{1}yu\partial_{u}),
\\
T_{1}=e^{-\gamma_{1}kt}(\partial_{t}+\gamma_{1}u\partial_{u}),\quad  T_{2}=e^{-\gamma_{1}
t}(\partial_{t}+\gamma_{1}\partial_{u}),\quad
\\
{\cal G}_{1}=e^{\gamma_{1}t}(\partial_{x}-\gamma_{1}\partial_{u}),\quad {\cal G}_{2}=e^{\gamma_{1}t}(\partial_{x}-\gamma_{1}u\partial_{u}),
\\
G_{0}=t\partial_{x}-\partial_{u},\quad
G_{1}=t\partial_{x}-u\partial_{u},
\\
Y=e^{t-\gamma_{1}x}u\partial_{u}
$,
\\
$R_{1}=e^{-x}(\cos y\partial_{x}-\sin y\partial_{y}
-\frac{2}{k}\cos y
u\partial_{u}),
\\
R_{2}=e^{-x}(\sin y\partial_{x}+\cos y\partial_{y}
-\frac{2}{k}\sin y
u\partial_{u})$,
\\
$Q^{1}_{\infty}=b(t,x,y)\partial_{u},\quad b(t,x,y)$  is an arbitrary solution
of the linear heat equation $b_{t}=\triangle b$,
\\
$Q^{2}_{\infty}=\beta(t,x,y)\partial_{u},\quad
\beta(t,x,y)$  is an  arbitrary solution
of the linear heat equation with the linear  source
$\beta_{t}=\triangle \beta \pm \beta$, \\
$X_{\infty}=A(x,y)\partial_{x}+B(x,y)\partial_{y}-2uA_{x}\partial_{u}$,
where $A(x,y)$  and $B(x,y)$  are arbitrary  functions satisfying the Cauchy-Riemann system $A_{x}=B_{y},\quad
A_{y}=-B_{x}$.
\end{remark}

\textbf{Proof.} First of all, we note that Cases 1 and 2 present
two subclasses of the general class of RDC equations (\ref{l1}),
when $K^{1}=K^{2}=0$. LSC of the first subclass was firstly derived in \cite{do-kn-sv-83}, while the second subclass was examined earlier in \cite{na-70}. So, Cases 3--16 present the results derived earlier in \cite{na-70}  and \cite{do-kn-sv-83}. They can be formally identified by examination of Case 1 from Table \ref{Tabl1-1}.

Examination of Cases 3,6,8 and 9 of Table \ref{Tabl1-1} is rather simple because $K^2=0$ and the function $K^1$ possesses a simple structure in each case. Notably these cases lead to the results derived in \cite{ed-broa-94} if one assumes  additionally that the relevant lambda-s vanish. The detailed  analysis involving all possible values of  lambda-s lead to the results listed in Cases 18, 22, and 28--32
of Table \ref{Tabl2-2}.

The equations arising in Cases 2,4,5 and 7 of Table \ref{Tabl1-1}  are absolutely new and their examination is very nontrivial.
Here we present the  detailed examination of Case 7.
In this case, the coefficients of Eq.(\ref{l1}) are specified as follows
\begin{equation}\begin{array}{l}\label{du20}
D=u^{k},\quad K^{1}=u^{k}\cos(p\ln u),\quad
K^{2}=u^{k}\sin(p\ln u),\quad
R=u(\lambda_{3}u^{k}+\lambda_{4}), k\neq 0.
\end{array}\end{equation}
Substituting the functions from
 (\ref{du20}) into the subsystem of DEs (\ref{m5})-(\ref{m7}), one obtains
\begin{equation}\begin{array}{l}\label{du21}
b=0,\quad ka=2\xi^{1}_{x}-\xi^{0}_{t},\quad \xi^{1}_{t}=0,\quad \xi^{2}_{t}=0,
\end{array}\end{equation}
\begin{equation}\begin{array}{l}\label{du22}
\xi^{1}_{x}\cos(p\ln u)-(pa-\xi^{2}_{x})\sin(p\ln u)=-2(k+1)a_{x},
\end{array}\end{equation}
\begin{equation}\begin{array}{l}\label{du2212}
(pa-\xi^{2}_{x})\cos(p\ln u)+\xi^{1}_{x}\sin(p\ln u)=-2(k+1)a_{y},
\end{array}\end{equation}
\begin{equation}\begin{array}{l}\label{du23}
a_{x}\cos(p\ln u)+a_{y}\sin(p\ln u)+2\lambda_{3}\xi^{1}_{x}=0,
\end{array}\end{equation}
\begin{equation}\begin{array}{l}\label{du24}
a_{t}=\lambda_{4}\xi^{0}_{t}.
\end{array}\end{equation}
Obviously, the general solution of  (\ref{du21})-(\ref{du24}) essentially depends on the parameter  $p$, hence two subcases {\it (i)} $p\neq0$ and {\it (ii)} $p=0$ should be examined.

Assuming
 $p\neq0$  and solving Eqs.(\ref{du21})-(\ref{du24}), we arrive at the linear system of the first-order PDE
\begin{equation}\begin{array}{l}\label{du25}
b=0,\quad a_{x}=0,\quad a_{y}=0,\quad  \xi^{1}_{t}=0,\quad \xi^{2}_{t}=0,\quad \xi^{1}_{x}=0,\quad  \xi^{2}_{x}=pa, \quad ka+\xi^{0}_{t}=0,
\end{array}\end{equation}
\begin{equation}\begin{array}{l}\label{du26}
\lambda_{4}\xi^{0}_{t}=0.
\end{array}\end{equation}
Now we consider two possible possibilities.

If $\lambda_{4}\neq0$ then Eq.(\ref{du26}) gives $\xi^{0}_{t}=0$, hence Eqs.(\ref{du25}) produce
$a=b=0,\quad \xi^{1}_{t}=\xi^{2}_{t}=\xi^{1}_{y}=\xi^{2}_{x}=0$.
 As a result, we realize that this possibility leads only to the principal algebra  $A^{pr}$.

If $\lambda_{4}=0$ then the general solution of system (\ref{du25})-- (\ref{du26}) (note that  Eqs.
(\ref{m3})-- (\ref{m4}) should be also taking into account) has the form

\begin{equation}\begin{array}{l}\label{du27}
\xi^{0}=kc_{0}t+d_{0},\quad \xi^{1}=pc_{0}y+d_{a},\quad \xi^{2}=-pc_{0}x+d_{a},\quad  \eta=-c_{0}u.
\end{array}\end{equation}

 The operator $X$ (\ref{m2})  with coefficients (\ref{du27}) produces the Lie algebra with the basic operators
\[\langle \partial_{t}, \partial_{x}, \partial_{y}, D_{3}+pJ_{12}\rangle.\]
So, Case 23 of Table \ref{Tabl2-2} is  identified  (an arbitrary $\lambda_{3}$  can be reduced to the three values $\sigma$ by a ET from $\mathcal{E}$).

Consider subcase {\it (ii)} $p=0$. Now
Eq. (\ref{l1}) with  (\ref{du20}) takes the form
\begin{equation}\begin{array}{l}\label{du2714}
u_{t}=(u^{k}u_{x})_{x}+(u^{k}u_{y})_{y}+u^{k}u_{x}+u(\lambda_{3}u^{k}+\lambda_{4}).
\end{array}\end{equation}

In order to simplify further calculations, we apply the following ET
of the form  (\ref{k1})
\begin{equation}\begin{array}{l}\label{du2715}
\frac{k^{2}}{16(k+1)^{2}}t\longrightarrow
t,\;
\frac{k}{4(k+1)}x\longrightarrow x,\;
\frac{k}{4(k+1)}y\longrightarrow
y,\;  u\longrightarrow u,
\\
16\frac{(k+1)^{2}}{k^{2}}\lambda_{3}\longrightarrow\lambda_{3},\; 16\frac{(k+1)^{2}}{k^{2}}\lambda_{4}\longrightarrow\lambda_{4},
\end{array}\end{equation}
which transform Eq. (\ref{du2714}) into
\begin{equation}\begin{array}{l}\label{du2716}
u_{t}=(u^{k}u_{x})_{x}+(u^{k}u_{y})_{y}+4\frac{k+1}{k}u^{k}u_{x}+u(\lambda_{3}u^{k}+\lambda_{4}),\quad k\neq -1; 0.
\end{array}\end{equation}

Thus, taking  into  account formulae  (\ref{du2715}),  equations  (\ref{du21})-(\ref{du24}) with  $p=0$ reduce to the form
\begin{equation}\begin{array}{l}\label{du28}
b=\xi^{1}_{t}=\xi^{2}_{t}=0,
\end{array}\end{equation}
\begin{equation}\begin{array}{l}\label{du29}
a=\frac{1}{2(k+1)}\xi^{1}_{x}-\frac{k}{16(k+1)^2}\xi^{0}_{t},
\end{array}\end{equation}
\begin{equation}\begin{array}{l}\label{du30}
2(k+1)a_{x}=-\xi^{1}_{x},
\end{array}\end{equation}
\begin{equation}\begin{array}{l}\label{du3097}
2(k+1)a_{y}=-\xi^{1}_{y},
\end{array}\end{equation}
\begin{equation}\begin{array}{l}\label{du31}
(\lambda_{3}-4\frac{k+1}{k^{2}})\xi^{1}_{x}=0,
\end{array}\end{equation}
\begin{equation}\begin{array}{l}\label{du32}
\xi^{0}_{tt}+k\lambda_{4}\xi^{0}_{t}=0.
\end{array}\end{equation}

The general solution of  (\ref{du28})-(\ref{du32}) essentially depends on the parameters  $\lambda_{3}$ and $ \lambda_{4}$. All the inequivalent subcases are

1)$\lambda_{3}\neq 4\frac{k+1}{k^{2}},\; \lambda_{4}=0$;

2)$\lambda_{3}\neq 4\frac{k+1}{k^{2}},\; \lambda_{4}\neq0$;

3)$\lambda_{3}=4\frac{k+1}{k^{2}},\; \lambda_{4}=0$;

4)$\lambda_{3}=4\frac{k+1}{k^{2}},\; \lambda_{4}\neq0$.

Consider  subcase 1)$\lambda_{3}\neq 4\frac{k+1}{k^{2}},\; \lambda_{4}=0$.
The general solution of system (\ref{du28})-(\ref{du32}) and  the remaining equations
(\ref{m3})-- (\ref{m4}) from the system of  DEs
is formed by the functions
\begin{equation}\begin{array}{l}\label{du33}
\xi^{0}=kc_{0}t+d_{0},\quad \xi^{a}=d_{a},\quad \eta=-c_{0}u.
\end{array}\end{equation}
 The operator $X$ (\ref{m2})with coefficients (\ref{du33}) produces the Lie algebra with the basic operators
\[\langle \partial_{t}, \partial_{x}, \partial_{y}, D_{3}\rangle.\]
So, Case 24 of Table \ref{Tabl2-2} is  identified.

Consider  subcase 2)$\lambda_{3}\neq 4\frac{k+1}{k^{2}}, \lambda_{4}\neq0$.
Using ET
\begin{gather*}
t\rightarrow\frac{\gamma_{1}}{\lambda_{4}}t,\quad
x\rightarrow x,\quad
y\rightarrow y,\quad
u\rightarrow \Big(\frac{\lambda_{4}}{\gamma_1}\Big)^{\frac{1}{k}}u.
\end{gather*}
we can make  $\lambda_{4}=\gamma_{1}$.
 So, the general solution of system (\ref{du28})-(\ref{du32}),  (\ref{m3}), (\ref{m4}) has the form
\begin{equation}\begin{array}{l}\label{du34}
\xi^{0}=c_{0}e^{-\gamma_{1}kt}+d_{0},\quad \xi^{a}=d_{a},\quad \eta=\gamma_{1}c_{0}e^{-\gamma_{1}kt}u.
\end{array}\end{equation}
 The operator $X$ (\ref{m2})with coefficients (\ref{du34}) produces the Lie algebra with the basic operators
\[\langle \partial_{t}, \partial_{x}, \partial_{y}, T_{1} \rangle.\]
So, Case 25 of Table \ref{Tabl2-2} is  identified.

Consider  subcase 3)$\lambda_{3}\neq 4\frac{k+1}{k^{2}},\; \lambda_{4}=0$.
 The linear ODE
\[ \xi^{1}_{xx}+\xi^{1}_{x}=0. \]
can be derived from Eqs.
(\ref{du29})-(\ref{du31}) in this subcase, which posses the general solution
\begin{equation}\begin{array}{l}\label{du36}
\xi^{1}=\varphi(y)e^{-x}+\psi(y).
\end{array}\end{equation}
Moreover, taking into account Eqs.
(\ref{m4})and  (\ref{du36})  we obtain
\begin{equation}\label{du38}
\xi^{2}=\dot{\varphi}(y)e^{-x}-\dot{\psi}(y)x+\chi(y),
\end{equation}
where  $\varphi=\varphi(y)$,   $ \psi=\psi(y)$and $\chi=\chi(y)$  are arbitrary functions at the moment.

Now we substitute  (\ref{du36})-(\ref{du38})  into Eqs. (\ref{m4}) and arrive at the linear ODE system
\begin{equation}\begin{array}{l}\label{du3812}
\varphi_{yy}+\varphi=0,\quad \psi_{yy}=0,\quad \chi_{yy}=0,
\end{array}\end{equation}
which possesses the general solution
\begin{gather*}
\varphi=c_{1}\cos y+c_{2}\sin y,\quad \psi=c_{3}y+c_{4},\quad \chi=c_{5}.
\end{gather*}
 So, we arrive at
\begin{equation}\begin{array}{l}\label{du39}
   \xi^{1}=e^{-x}(c_{1}\cos y+c_{2}\sin y)+c_{3}y+d_{1},
\end{array}\end{equation}
\begin{equation}\begin{array}{l}\label{du40}
   \xi^{2}=e^{-x}(-c_{1}\sin y+c_{2}\cos y)-c_{3}x+d_{2},
\end{array}\end{equation}

The coefficient $\xi^{0}$ of the operator  $X$ can be easily derived from Eq.
(\ref{du32}):
\begin{equation}\begin{array}{l}\label{du41}
\xi^{0}=kc_{0} t+d_{0}.
\end{array}\end{equation}
 In the above formulae $c_{0}, c_{1},\quad  c_{2}, \quad  c_{3},\quad  c_{4}, \quad  c_{5}, \quad d_{0}, \quad d_{1}$  and  $d_{2}$ are arbitrary constants.

In order to find the  coefficient $\xi^{0}$ of the operator  $X$, we use (\ref{du29}), (\ref{du41}):
\begin{equation}\begin{array}{l}\label{du4212}
a=-\frac{1}{2(k+1)}e^{-x}(c_{1}\cos y+c_{2}\sin y)-\frac{k}{16(k+1)^2}c_{0}.
\end{array}\end{equation}
So, taking into account  (\ref{3-1}), (\ref{du28})and  (\ref{du4212}),  we obtain
\begin{equation}\begin{array}{l}\label{du309}
\eta=\Big[\frac{1}{2(k+1)}e^{-x}(c_{1}\cos y+c_{2}\sin y) + \frac{k}{16(k+1)^2}c_{0}\Big]u.
\end{array}\end{equation}

 Substituting (\ref{du4212}) and  (\ref{du39}) into  (\ref{du3097}), we arrive at the restriction $c_{3}=0$.
So,
the operator $X$ (\ref{m2})with coefficients (\ref{du41}), (\ref{du39})--(\ref{du40}) (under the restriction$c_3=0$) and  (\ref{du309})  produces the Lie algebra with the basic operators
\[\langle \partial_{t}, \; \partial_{x}, \; \partial_{y},\;  D_{3},\; R_{1},\; R_{2} \rangle.\]
Thus, Case 26 of Table \ref{Tabl2-2} is  identified.

Finally, we examine   subcase 4)$\lambda_{3}\neq 4\frac{k+1}{k^{2}}, \lambda_{4}\neq0$ .
Using ET
\[t\longrightarrow\frac{t}{\lambda_{4}},\; x\longrightarrow x,\; y\longrightarrow y,\;  u\longrightarrow |\lambda_{4}|^{\frac{1}{k}}u,\]
we  can  set  $\lambda_{4}=\gamma_{1}$ without losing a generality.
 So, solving Eq. (\ref{du32}) we obtain
\begin{equation}\label{du43}
\xi^{0}=c_{0} e^{-\gamma_{1}k t}+d_{0}.
\end{equation}
Eqs.(\ref{3-1}), (\ref{du28}) and  (\ref{du43}) give
\begin{equation}\begin{array}{l}\label{du4312}
\eta=-\Big[\frac{1}{2(k+1)}e^{-x}(c_{1}\cos y+c_{2}\sin y)+\gamma_{1}c_{0}e^{-\gamma_{1}k t}\Big]u.
\end{array}\end{equation}
The coefficients $\xi^{1}$ and $\xi^{2}$ again are given by
(\ref{du39})--(\ref{du40}) with $c_3=0$.

The operator $X$ (\ref{m2})with coefficients (\ref{du43}), (\ref{du39})--(\ref{du40}) (under the restriction  $c_3=0$) and  (\ref{du4312})  produces the Lie algebra with the basic operators
\[\langle \partial_{t}, \; \partial_{x}, \; \partial_{y},\;  T_{1},\; R_{1},\; R_{2} \rangle.\]
So, Case 27 of Table  is  identified.

Thus, Cases 23-27   of Table \ref{Tabl2-2} have been   identified by examination of the  RDC equation  (\ref{l1}) with the coefficients  listed  in  Case 7 of Table \ref{Tabl1-1}.

  Cases 17-22 of Table \ref{Tabl2-2} have been   obtained by a similar  analysis of  the equations  with the coefficients  listed  in Cases 2-6 of Table \ref{Tabl1-1}.

  Finally,  Cases  28-32 of Table \ref{Tabl2-2} have been identified  by the analysis  of the RDC equations with the coefficients  listed   in Cases 8-9 of Table \ref{Tabl1-1}.

The proof is now completed.\quad $\Box$

Thus, we can state that the first five steps of the LSC algorithm presented in Section 2 have been realized. As a result, we have derived LSC of the class of RDC equations (\ref{l1}) based on the group of ETs (\ref{k1}). Such classification is often called LSC via the Lie-Ovsiannikov algorithm (see \cite{ch-se-pl-2018} for discussion  on this matter). However, it is well-known that the Lie-Ovsiannikov algorithm does not lead to the so called canonical list of the PDEs admitting nontrivial Lie symmetry. In fact, the number of relevant  equations often can be reduced by implementation of the last two steps of the algorithm from Section 2. In the next section, it will be proved that 10 equations among 32 those from Table \ref{Tabl2-2} are reducible to other equations from the same table by appropriate FPTs.

 \section {\bf Lie symmetry classification   using the form-preserving  transformations}

Now we turn to notion of a form-preserving transformation (FPT).
Roughly speaking,  a  FPT  is a local substitution, which reduces some PDE from the given class to another PDE belonging to the same class. The rigorous definition can be as follows (see \cite{ch-se-pl-2018}, P.32).

\textbf{Definition.}
A non-degenerate point transformation given by the formulae
\begin{equation}\label{2-4-1*}
{t}^* = f(t,x,u), \quad
x^*_a = g_a(t,x,u),\quad
u^* =
h(t,x,u),\ (a=1,\dots,n),
 \end{equation}
 which maps at least one equation of the form
(\ref{1-1-0}) into an equation belonging to the same class, is called
the FPT of the PDE class (\ref{1-1-0}).

\medskip

Comparing this definition with the well-known definition of ETs, one immediately notes that  each ET
is automatically a FPT but not vice versa.
In contrast to the ETs, a set of all possible
FPTs for the given class of PDEs usually do not form a Lie group. However, a subset of FPTs may generate a  group of ETs on a subclass of the given class (see, e.g. example in \cite{ch-se-pl-2018}, Section 2.3.2 ). This is a reason  why FPTs  are also called additional equivalence transformations. To the best of our knowledge, the 1992  paper  \cite{ga-wi-92}  was the first, in which  FPTs were used to solve LSC problem for a class  PDEs (the authors used the terminology `admissible transformations').

Let us construct the set of FPTs for the class of RDC equations (\ref{l1}). We start from the most general form of point transformations
\begin{equation}\begin{array}{l}\label{d1}
\tau=a(t,x,y,u),\quad x^{*}=b^{1}(t,x,y,u),\quad y^{*}=b^{2}(t,x,y,u),\quad v=c(t,x,y,u).
\end{array}\end{equation}
Now we  assume that there exists a FPT of the form (\ref{d1}), which relates an equation from the class (\ref{l1}) with another one from the same class, say, of the form
\begin{equation}\begin{array}{l}\label{d2}
v_{\tau}=(d(v)v_{x^{*}})_{x^{*}}+(d(v)v_{y^{*}})_{y^{*}}+k^{1}(v)v_{x^{*}}+k^{2}(v)v_{y^{*}}+r(v),
\end{array}\end{equation}
Here  $u=u(t,\;  x,\;  y)$ and $  v=v(\tau,\;
x^{*},\; y^{*})$ are unknown functions, while $a(t,x,y,u),\quad b^{1}(t,x,y,u)$, $ b^{2}(t,x,y,u),c(t,x,y,u)\; d=d(v),\;  k^{1}=k^{1}(v),\; k^{2}=k^{2}(v),\; r=r(v)$  are some given functions.

\begin{theo}\label{T2-5-5}
An arbitrary  RDC equation of the form   (\ref{l1}) can be reduced to another equation of the same form (\ref{d2}) by the local nondegenerate  transformation (\ref{d1}) with the correctly-specified smooth functions $a, b^1, b^2$ and $c$ if and only if these functions are of the form
\begin{equation}\begin{array}{l}\label{d3}
\tau=a(t),\; x^{*}=b^{1}(t,x,y),\; y^{*}=b^{2}(t,x,y),\;  v=M(t,x,y)u+N(t,x,y),
\end{array}\end{equation}
and and the following equalities take place
\begin{equation}\begin{array}{l}\label{d5}
b^{2}_{x}=\pm b^{1}_{y},
\end{array}\end{equation}
\begin{equation}\begin{array}{l}\label{d6}
b^{2}_{y}=\mp b^{1}_{x},
\end{array}\end{equation}
\begin{equation}\begin{array}{l}\label{d7}
\Big[(b^{1}_{x})^{2}+(b^{1}_{y})^{2}\Big]D(u)=\dot{a}d(v),
\end{array}\end{equation}
\begin{equation}\begin{array}{l}\label{d8}
b^{1}_{t}+\frac{2}{M}\frac{d}{du}\Big[\Big[b^{1}_{x}(M_{x}u+N_{x})+b^{1}_{y}(M_{y}u+N_{y})\Big]{D}(u)\Big]-b^{1}_{x}K^{1}(u)-b^{1}_{y}K^{2}(u)=-\dot{a}k^{1}(v),
\end{array}\end{equation}
\begin{equation}\begin{array}{l}\label{d814}
b^{2}_{t}+\frac{2}{M}\frac{d}{du}\Big[\Big[b^{2}_{x}(M_{x}u+N_{x})+b^{2}_{y}(M_{y}u+N_{y})\Big]D(u)\Big]-b^{2}_{x}K^{1}(u)+b^{2}_{y}K^{2}(u)=-\dot{a}k^{2}(v),
\end{array}\end{equation}
\begin{equation}\begin{array}{l}\label{d9}
M_{t}u+N_{t}-(\triangle M u+\triangle N)D(u)+\frac{1}{M}\frac{d}{du}\Big[\Big[(M_{x}u+N_{x})^2+(M_{y}u+N_{y})^2\Big]D(u)\Big]-
\\
-(M_{x}u+N_{x})K^{1}(u)-(M_{y}u+N_{y})K^{2}(u)+M
R(u)=\dot{a}r(v).
\end{array}\end{equation}
provided
\begin{equation}\begin{array}{l}\label{d4}
\dot{a}M\Big[(b^{1}_{x})^{2}+(b^{1}_{y})^{2}\Big]\neq 0.
\end{array}\end{equation}
\end{theo}

\textbf{Proof.}
Firstly we note that any FPT (\ref{d3}) must be nondegenerate, i.e., its Jacobian is nonvanish:
\begin{equation}\label{d10} J=\begin{vmatrix}
a_{t} &  a_{x} & a_{y} & a_{u}\\
b^{1}_{t} &  b^{1}_{x} & b^{1}_{y} & b^{1}_{u}\\
b^{2}_{t} &  b^{2}_{x} & b^{2}_{y} & b^{2}_{u}\\
c_{t} &  c_{x} & c_{y} & c_{u}
\end{vmatrix} \neq 0.
\end{equation}

Having transformation  (\ref{d1}) one can express the derivatives of the function $u$
by the well-known formulas (usually they are presented in the case of two independent variables but those formulae can be directly extended on three or more variables)

\begin{equation}\begin{array}{l}\label{d12}
u_{t}=-\frac{v_{\tau}a_{t}+v_{x^{*}} b^{1}_{t}+v_{y^{*}} b^{2}_{t}-
c_{t}}{A},
\end{array}\end{equation}
\begin{equation}\begin{array}{l}\label{d130}
u_{x}=-\frac{v_{\tau}a_{x}+v_{x^{*}} b^{1}_{x}+v_{y^{*}} b^{2}_{x} -
c_{x}}{A},
\end{array}\end{equation}
\begin{equation}\begin{array}{l}\label{d132}
u_{y}=-\frac{v_{\tau}a_{y}+v_{x^{*}} b^{1}_{y}+v_{y^{*}} b^{2}_{y} -
c_{y}}{A},
\end{array}\end{equation}
\begin{equation}\begin{array}{l}\label{d140}
u_{xx}=-\frac{1}{A}\Big[(a_{x}+a_{u}u_{x})^{2}v_{\tau \tau}+2(a_{x}+a_{u}u_{x})(b^{1}_{x}
+b^{1}_{u}u_{x})v_{\tau x^{*}}
+2(a_{x}+a_{u}u_{x})(b^{2}_{x}+
\\ [1mm]
+b^{2}_{u}u_{x})v_{\tau y^{*}}+
(b^{1}_{x}+b^{1}_{u}u_{x})^{2}v_{x^{*}}v_{x^{*}}+(b^{1}_{x}+b^{1}_{u}u_{x})(b^{2}_{x}+b^{2}_{u}u_{x})v_{x^{*}}v_{y^{*}}+
(b^{2}_{x}+b^{2}_{u}u_{x})^{2}v_{y^{*}}v_{y^{*}}+
\\ [1mm]
+(a_{xx}+2a_{xu}u_{x}
+a_{uu}u_{x}u_{x})v_{\tau}
+(b^{1}_{xx}+2b^{1}_{xu}u_{x}+
b^{1}_{uu}u_{x}u_{x})v_{x^{*}}
+(b^{2}_{xx}+2b^{2}_{xu}u_{x}+
\\ [1mm]
+b^{2}_{uu}u_{x}u_{x})v_{y^{*}}-
(c_{xx}+2c_{x
u}u_{x}+c_{uu}u_{x}u_{x})\Big],
\end{array}\end{equation}
\begin{equation}\begin{array}{l}\label{d141}
u_{yy}=-\frac{1}{A}\Big[(a_{y}+a_{u}u_{y})^{2}v_{\tau \tau}+2(a_{y}+a_{u}u_{y})(b^{1}_{y}
+b^{1}_{u}u_{y})v_{\tau x^{*}}
+2(a_{y}+a_{u}u_{y})(b^{2}_{y}+
\\ [1mm]
+b^{2}_{u}u_{y})v_{\tau y^{*}}+
(b^{1}_{y}+b^{1}_{u}u_{y})^{2}v_{x^{*}}v_{x^{*}}+(b^{1}_{y}+b^{1}_{u}u_{x})(b^{2}_{y}+b^{2}_{u}u_{y})v_{x^{*}}v_{y^{*}}+
(b^{2}_{y}+b^{2}_{u}u_{y})^{2}v_{y^{*}}v_{y^{*}}+
\\ [1mm]
+(a_{yy}+2a_{yu}u_{y}
+a_{uu}u_{y}u_{y})v_{\tau}
+(b^{1}_{yy}+2b^{1}_{yu}u_{x}+
b^{1}_{uu}u_{y}u_{y})v_{x^{*}}
+(b^{2}_{yy}+2b^{2}_{yu}u_{y}+
\\ [1mm]
+b^{2}_{uu}u_{y}u_{y})v_{y^{*}}-
(c_{yy}+2c_{y
u}u_{y}+c_{uu}u_{y}u_{y})\Big],
\end{array}\end{equation}
where $ A=v_{\tau}a_{u}+v_{x^{*}} b^{1}_{u}+v_{y^{*}} b^{2}_{u}- c_{u},\;
u_{x}=\frac{\partial u}{\partial x},\; u_{y}=\frac{\partial u}{\partial y},\;
v_{x^{*}}=\frac{\partial v}{\partial x^{*}},\quad v_{y^{*}}=\frac{\partial v}{\partial y^{*}},\quad
u_{xx}=\frac{\partial^{2}u}{\partial x^{2}},\; u_{yy}=\frac{\partial^{2}u}{\partial y^{2}},\;
v_{x^{*}x^{*}}=\frac{\partial^{2}v}{\partial (x^{*})^{2}}, v_{x^{*}y^{*}}=\frac{\partial^{2}v}{\partial x^{*} \partial y^{*}}, v_{y^{*}y^{*}}=\frac{\partial^{2}v}{\partial (y^{*})^{2}}$.

Substituting (\ref{d12})--(\ref{d141}) into (\ref{l1}) one arrives at a very cumbersome expression.
Let us assume that (\ref{d1}) is a FPT.
So,    the expression obtained must be  reducible to an equation of the form (\ref{d2}).
In the particular case, the coefficient next to the second-order derivative $v_{\tau \tau},\; v_{\tau x^{*}},\; v_{\tau y^{*}}$ and  $v_{x^{*} y^{*}}$
should vanish and  the coefficients next to the  derivatives $v_{x^{*} x^{*}}$  and $ v_{y^{*} y^{*}}$ must be equal to $u_{xx}$ and $v_{xx}$, respectively, hence
one obtains the system of PDEs

\begin{equation}\begin{array}{l}\label{d15}
a_{x}+a_{u}u_{x}=0,\; a_{y}+a_{u}u_{y}=0,
\end{array}\end{equation}
\begin{equation}\begin{array}{l}\label{d16}
(b^{1}_{x}+b^{1}_{u}u_{x})(b^{2}_{x}+b^{2}_{u}u_{x})=0, \; (b^{1}_{y}+b^{1}_{u}u_{y})(b^{2}_{y}+b^{2}_{u}u_{y})=0,
\end{array}\end{equation}
\begin{equation}\begin{array}{l}\label{d17}
(b^{1}_{x}+b^{1}_{u}u_{x})(b^{1}_{x}+b^{1}_{u}u_{x})=(b^{2}_{x}+b^{2}_{u}u_{x})(b^{2}_{x}+b^{2}_{u}u_{x}),\;
\\
(b^{1}_{y}+b^{1}_{u}u_{y})(b^{1}_{y}+b^{1}_{u}u_{y})=(b^{2}_{y}+b^{2}_{u}u_{y})(b^{2}_{y}+b^{2}_{u}u_{y}).
\end{array}\end{equation}

Equations  (\ref{d15}) immediately give
\begin{equation}\begin{array}{l}\label{d18}
a_{x}=a_{y}=a_{u}=0,
\end{array}\end{equation}
while Eqs. (\ref{d17})  lead to
\begin{equation}\begin{array}{l}\label{d19}
b^{1}_{x}b^{2}_{x}+b^{1}_{y}b^{2}_{y}=0,
\end{array}\end{equation}
\begin{equation}\begin{array}{l}\label{d20}
b^{1}_{x}b^{2}_{u}+b^{1}_{u}b^{2}_{x}=0,\; b^{1}_{y}b^{2}_{u}+b^{1}_{u}b^{2}_{y}=0,
\end{array}\end{equation}
\begin{equation}\begin{array}{l}\label{d21}
b^{1}_{u}b^{2}_{u}=0,
\end{array}\end{equation}
\begin{equation}\begin{array}{l}\label{d22}
(b^{1}_{x})^{2}+(b^{1}_{y})^{2}=(b^{2}_{x})^{2}+(b^{2}_{y})^{2},
\end{array}\end{equation}
\begin{equation}\begin{array}{l}\label{d23}
b^{1}_{x}b^{1}_{u}=b^{2}_{x}b^{2}_{u},\; b^{1}_{y}b^{1}_{u}=b^{2}_{y}b^{2}_{u},
\end{array}\end{equation}
\begin{equation}\begin{array}{l}\label{d24}
(b^{1}_{u})^{2}=(b^{2}_{u})^{2}.
\end{array}\end{equation}
Obviously, Eqs. (\ref{d21})and  (\ref{d24}) are equivalent to
\begin{equation}\begin{array}{l}\label{d25}
b^{1}_{u}=b^{2}_{u}=0.
\end{array}\end{equation}
So, using  the derived restrictions (\ref{d18}) and (\ref{d25}) we can specify FPT in question as follows
\begin{equation}\begin{array}{l}\label{d26}
\tau=a(t),\; x^{*}=b^{1}(x,y),\; y^{*}=b^{2}(x,y),\; v=c(x,y,u),
\end{array}\end{equation}
i.e. the first three formulae in (\ref{d2}) are derived.

Now one observes that formulae
(\ref{d12})-(\ref{d141})  can be simplified essentially if one takes into account  (\ref{d25}), namely:
\begin{equation}\begin{array}{l}\label{d27}
u_{t}=\frac{1}{c_{u}}(\dot{a}v_{\tau}+b^{1}_{t}v_{x^{*}}+b^{2}_{t}v_{y^{*}}-c_{t}),
\end{array}\end{equation}
\begin{equation}\begin{array}{l}\label{d2811}
u_{x}=\frac{1}{c_{u}}(b^{1}_{x}v_{x^{*}}+b^{2}_{x}v_{y^{*}}-c_{x}),
\end{array}\end{equation}
\begin{equation}\begin{array}{l}\label{d2812}
u_{y}=\frac{1}{c_{u}}(b^{1}_{y}v_{x^{*}}+b^{2}_{y}v_{y^{*}}-c_{y}),
\end{array}\end{equation}
\begin{equation}\begin{array}{l}\label{d29}
\triangle u=\frac{1}{c_{u}}\Big[\Big((b^{1}_{x})^{2}+(b^{1}_{y})^{2}\Big)\triangle
v-\frac{c_{uu}}{c^{2}_{u}}\Big[\Big((b^{1}_{x})^{2}+(b^{1}_{y})^{2}\Big)v_{x^{*}}v_{x^{*}}+2(b^{1}_{x}b^{2}_{x}+b^{1}_{y}b^{2}_{y})v_{x^{*}}v_{y^{*}}+\Big((b^{2}_{x})^{2}+
\\
+(b^{2}_{y})^{2}\Big)v_{y^{*}}v_{y^{*}}\Big]
+\Big(\triangle
b^{1}-\frac{2c_{xu}}{c_{u}}b^{1}_{x}-\frac{2c_{yu}}{c_{u}}b^{1}_{y}+2\frac{c_{uu}}{c^{2}_{u}}(b^{1}_{x}c_{x}+b^{1}_{y}c_{y})\Big)v_{x^{*}}
+\Big(\triangle
b^{2}-\frac{2c_{xu}}{c_{u}}b^{2}_{x}-
\\
-\frac{2c_{yu}}{c_{u}}b^{2}_{y}+2\frac{c_{uu}}{c^{2}_{u}}(b^{2}_{x}c_{x}+b^{2}_{y}c_{y})\Big)v_{y^{*}}
+2\frac{c_{xu}}{c_{u}}c_{x}+2\frac{c_{yu}}{c_{u}}c_{y}-\frac{c_{uu}}{c^{2}_{u}}(c_{x}c_{x}+c_{x}c_{y}+c_{y}c_{y})-\triangle
c\Big].
\end{array}\end{equation}

Finally, substituting the right-hand-sides from  (\ref{d27})-(\ref{d29}) into  (\ref{l1}), we note that the expression obtained is reducible to Eq. (\ref{d2}) only under the condition
\begin{equation}\begin{array}{l}\label{d31}
c_{uu}=0,
\end{array}\end{equation}
i.e. the last formula in (\ref{d2}) is derived,
and the equalities (\ref{d7})-(\ref{d9}) should take place.

It can be also easily shown that Eqs.(\ref{d19}) and (\ref{d22}) are equivalent to Eqs.(\ref{d5})--(\ref{d6}), while the restriction (\ref{d4}) immediately follows from (\ref{d10}) because of formulae (\ref{d3}).

The proof is now completed.\quad $\Box$


\begin{remark}
Using Theorem \ref{T2-5-5}, one can derive the discrete equivalence transformation

\begin{equation}\begin{array}{l}\label{d122}
t^{*}=t,\; x^{*}=-x,\; y^{*}=y,\; u^{*}=u;
\end{array}\end{equation}
\begin{equation}\begin{array}{l}\label{d123}
t^{*}=t,\; x^{*}=x,\; y^{*}=-y,\; u^{*}=u;
\end{array}\end{equation}
\begin{equation}\begin{array}{l}\label{d121}
t^{*}=-t,\; x^{*}=x,\; y^{*}=y,\; u^{*}=u;
\end{array}\end{equation}
\begin{equation}\begin{array}{l}\label{d124}
t^{*}=t,\; x^{*}=x,\; y^{*}=y,\; u^{*}=-u;
\end{array}\end{equation}
\begin{equation}\begin{array}{l}\label{d1248}
t^{*}=t,\; x^{*}=y,\; y^{*}=x,\; u^{*}=u,
\end{array}\end{equation}
which  take place for  Eq. (\ref{l1}) with arbitrary smooth coefficients     $D, \ K^{1},\; K^{2},\; R$.  However, if one takes into account the  restriction   $D(u)>0$ then the discrete  transformation (\ref{d121}) is not valid, while (\ref{d124}) is valid only under the additional  restriction $D(-u)>0$.
\end{remark}

To complete the last step of the  algorithm (see Section 2), we need to apply Theorem \ref{T2-5-5} to the equations listed in Table \ref{Tabl2-2} in order to identify those pairs of them, which are reducible one to another by an appropriate FPT.

\begin{theo}\label{T2-5-6} There are
9 equations in  Table \ref{Tabl2-2}, which are reducible to other equations from the same table
by an appropriate FPT of the form \eqref{d3}. All the equations and the corresponding transformations are presented in Cases 1--9 of Table \ref{Tabl3-1}.
\end{theo}

\begin{longtable}[�]{|c|c|c|c|}
\caption[�]{Simplification of the RDC equations form (\ref{l1})  by means  of FPTs}\label{Tabl3-1}\\
\hline  & RDC equation & FPT & Canonical form
\\ [0.5mm]
&  & & of RDC equation  \\ [0.5mm]
\hline 1 & $u_{t}=\triangle u+\gamma_{1}$&$\tau=t,\; x^{*}=x,$& $v_{\tau}=\triangle v$
\\ [0.5mm]
& $$&$y^{*}=y,\; v=u
-\gamma_{1}t$& $$\\ [0.5mm]
\hline 2 & $u_{t}=\triangle u+\gamma_{1}u$&$\tau=t,\; x^{*}=x,$& $v_{\tau}=\triangle v$
\\ [0.5mm]
& $$&$y^{*}=y,\;
v=e^{\gamma_{1}t}u$& $$
\\ [0.5mm]
\hline 3 & $u_{t}=(e^{u}u_{x})_{x}+(e^{u}u_{y})_{y}+$ & $\tau=\gamma_{1}e^{\gamma_{1}
t},\; x^{*}=x,$& $v_{\tau}=(e^{v}v_{x^{*}})_{x^{*}}+(e^{v}v_{y^{*}})_{y^{*}}+\sigma e^{v}$
\\ [0.5mm]
&$+\sigma e^{u}+\gamma_{1}$&$y^{*}=y,\; v=u+\gamma_{1}t$&$$
\\ [0.5mm]
\hline 4 & $u_{t}=(u^{k}u_{x})_{x}+(u^{k}u_{y})_{y}+$&$\tau=\gamma_{1}\frac{e^{\gamma_{1}kt}}{k},\;
x^{*}=x,$& $v_{\tau}=(v^{k}u_{x^{*}})_{x^{*}}+(v^{k}u_{y^{*}})_{y^{*}}+\sigma v^{k+1}$
\\ [0.5mm]
&$+\sigma u^{k+1}+\gamma_{1}u$&$y^{*}=y,\; v=e^{-\gamma_{1}t}u$&$$
\\ [0.5mm]
\hline 5 & $u_{t}=(e^{u}u_{x})_{x}+(e^{u}u_{y})_{y}+$ &
$\tau=\gamma_{1}e^{\gamma_{1}t},\; x^{*}=x,$& $v_{\tau}=(e^{v}u_{x^{*}})_{x^{*}}+(e^{v}u_{y^{*}})_{y^{*}}+$
\\ [1mm]
&$+e^{u}u_{x}+\sigma e^{u}+\gamma_{1}$
&$y^{*}=y,\; v=u-\gamma_{1}t$&$+e^{v}v_{x^{*}}+\sigma e^{v}$
\\ [1mm]
\hline 6 & $u_{t}=(u^{k}u_{x})_{x}+(u^{k}u_{y})_{y}+$&$\tau=\gamma_{1}\frac{e^{\gamma_{1}kt}}{k},\;
x^{*}=x,$& $v_{\tau}=(v^{k}u_{x^{*}})_{x^{*}}+(v^{k}u_{y^{*}})_{y^{*}}+$
\\ [1mm]
&$+u^{k}u_{x}+\sigma u^{k+1}+\gamma_{1}u$&$y^{*}=y,\;  v=e^{-\gamma_{1}t}u$
&$+v^{k}v_{x^{*}}+\sigma v^{k+1}$
\\ [1mm]
\hline 7 & $u_{t}=(u^{k}u_{x})_{x}+(u^{k}u_{y})_{y}+$&
$\tau=\gamma_{1}\frac{e^{\gamma_{1}kt}}{k}, \; x^{*}=x,$& $v_{\tau}=(v^{k}u_{x^{*}})_{x^{*}}+(v^{k}u_{y^{*}})_{y^{*}}+$
\\ [1.5mm]
&
$+4\frac{k+1}{k}u^{k}u_{x}+4\frac{k+1}{k^{2}}u^{k+1}+\gamma_{1}u$&$y^{*}=y,\; v=e^{-\gamma_{1}
t}u$
&
$+4\frac{k+1}{k}v^{k}v_{x^{*}}+4\frac{k+1}{k^{2}}v^{k+1}$\\ [1.5mm]
\hline 8 & $u_{t}=\triangle u+uu_{x}+\sigma$ &
$\tau=t, \; x^{*}=x+\gamma_{1}\frac{1}{2}t^{2},$& $v_{\tau}=\triangle
v+vv_{x^{*}}$
\\ [1mm]
&$$&$y^{*}=y,\; v=u-\gamma_{1}
t$&$$
\\ [1mm]
\hline 9 & $u_{t}=\triangle u+\ln uu_{x}+\sigma u$ &
$\tau=t, \; x^{*}=x+\gamma_{1}\frac{1}{2}t^{2},$& $v_{\tau}=\triangle
v+\ln vv_{x^{*}}$
\\ [1mm]
&$$&$y^{*}=y,\; v=e^{-\gamma_{1}t}u$&$$
\\ [1mm]
\hline
\end{longtable}

\textbf{Sketch of the proof} of this theorem is similar to  that of Theorem 2.11 \cite{ch-se-pl-2018}. Of course, Cases 1--4 of Table \ref{Tabl3-1} involve the well-known substitutions for the linear and nonlinear RD, which were identified  many years ago.

The peculiarity of the RDC equations  in Cases 5--9 is such that each equation contains  the convective, which  involves only the derivative w.r.t. $x$  (no terms involving $u_y$). It turns out that FPTs constructed for the (1+1)-dimensional analogs of these equations in \cite{ch-se-2006} (see also Table 2.6 in \cite{ch-se-pl-2018}) are valid also for (1+2)-dimensional equations (the second space variable $y$ is unchangeable).

In particular, the most nontrivial FPT occurs in Case 9 of Table \ref{Tabl3-1}. On the other hand, one notes that it is nothing else but the substitution listed in Case 16 of Table 2.6  \cite{ch-se-pl-2018} with the formal additional transformation $y^{*}=y$.  Substituting these formulae into Eqs. (\ref{d5})-(\ref{d9}) and taking into account that
\[ D(u)=d(v)=1, \  K^1=\ln u, \  k^1=\ln v, \  K^2=k^2=0, \ R(u)=\sigma u, \ r(v)=0, \]
one easily checks that the substitution
\[\tau=t, \ x^{*}=x+\frac{1}{2}\sigma t^{2}, \ y^{*}=y,\ v=e^{-\sigma t}u \]
is indeed FPT, which  relates  two RDC equations listed   in Case 9 of Table \ref{Tabl3-1}.


The proof is now completed.\quad $\Box$

Now we formulate the main theorem  presenting the canonical list of RDC equations from class (\ref{l1}) possessing nontrivial Lie symmetries.

\begin{theo}\label{T2-5-7}
 All possible  RDC equations of the form (\ref{l1}) admitting  nontrivial Lie symmetries  are reduced to one of the 22 canonical  equations  listed  in the second column of Table \ref{Tabl4-1}
 by the relevant FPTs presented in Theorem \ref{T2-5-5}.
 The relevant  MAIs  of the canonical RDC equations
 are listed in the third column of Table  \ref{Tabl4-1}.
 \end{theo}

\begin{longtable}[c]{|c|c|c|c|c|}
\caption[c]{The full LSC  of the class of  RDC equations  (\ref{l1})}\label{Tabl4-1}\\
\hline
& Equation  & MAI & Restrictions  \\ [0.5mm]
\hline 1 & $u_{t}=(D(u)u_{x})_{x}+(D(u)u_{y})_{y}+$& $<\partial_{t},
\partial_{x}, \partial_{y}, J_{12}>$ & $D-\forall, R-\forall$ \\ [0.5mm]
&$+R(u)$& $$ & $$ \\ [0.5mm]
\hline 2 & $u_{t}=(D(u)u_{x})_{x}+(D(u)u_{y})_{y}$ & $<\partial_{t},
\partial_{x}, \partial_{y}, J_{12}, D_{0}>$ & $D-\forall$ \\ [0.5mm]
\hline 3 & $u_{t}=\triangle u$ & $<\partial_{t},
\partial_{x}, \partial_{y}, J_{12}, G_{x}, G_{y}, I,$ & $$ \\ [0.5mm]
&& $D_{0}, \Pi,
 Q^{1}_{\infty}>$ & $$
\\ [0.5mm]
\hline 4 & $u_{t}=\triangle u+\gamma_{1}u \ln u$ & $<\partial_{t},
\partial_{x}, \partial_{y},  J_{12}, e^{\gamma_{1} t}I, {\cal G}_{x}, {\cal G}_{y}>$ & $$
\\ [0.5mm]
\hline 5 & $u_{t}=(e^{u}u_{x})_{x}+(e^{u}u_{y})_{y}$ & $<\partial_{t},
\partial_{x}, \partial_{y}, J_{12}, D_{0}, D_{2}>$ & $\delta=1$
\\ [0.5mm]
\hline 6 & $u_{t}=(e^{\delta u}u_{x})_{x}+(e^{\delta u}u_{y})_{y}+\gamma_{1}e^{mu}$ & $<\partial_{t},
\partial_{x}, \partial_{y}, J_{12}, $ & $m \neq 0$ \\ [0.5mm]
&& $(\delta-m)D_{0}-2 D_{4}>$ & $$
\\ [0.5mm]
\hline 7 & $u_{t}=(u^{k}u_{x})_{x}+(u^{k}u_{y})_{y}$ & $<\partial_{t},
\partial_{x}, \partial_{y}, J_{12}, D_{0}, D_{1}>$ & $k \neq -1; 0$\\ [0.5mm]
\hline 8 & $u_{t}=(u^{k}u_{x})_{x}+(u^{k}u_{y})_{y}+\gamma_{1}u^{m}$ & $<\partial_{t},
\partial_{x}, \partial_{y}, J_{12}, $ & $k \neq -1; 0$ \\ [0.5mm]
&& $(m-1)D_{0}-D_{1}>$
&$m \neq 1$\\ [0.5mm]
\hline 9 & $u_{t}=(u^{-1}u_{x})_{x}+(u^{-1}u_{y})_{y}$ & $<\partial_{t},
\partial_{x}, \partial_{y}, J_{12}, D_{3},X_{\infty}>$ & $k=-1$
\\ [0.5mm]
\hline 10& $u_{t}=(e^{\delta
u}u_{x})_{x}+(e^{\delta
u}u_{y})_{y}+$&
$<\partial_{t}, \partial_{x},\partial_{y},$&$m\neq \delta$\\ [0.5mm]
&$+e^{mu}[u_{x}\cos(pu)+u_{y}\sin(pu)]+\sigma e^{(2m-\delta)u}$
&$(m-\delta)D_{0}+D_{4}+pJ_{12}>$&$(m,p)\neq (0,0)$\\ [0.5mm]
\hline 11 & $u_{t}=(e^{u}u_{x})_{x}+(e^{u}u_{y})_{y}+$ & $<\partial_{t}, \partial_{x},\partial_{y}, $& $\delta=1$\\ [0.5mm]
&$+uu_{x}+\sigma e^{-u}$&$D_{0}-D_{4}-t\partial_{x}>$&$$
\\ [0.5mm]
\hline 12 & $u_{t}=(e^{u}u_{x})_{x}+(e^{u}u_{y})_{y}+$ & $<\partial_{t}, \partial_{x}, \partial_{y},
D_{4}+pJ_{12}>$ & $\delta=1$
\\ [0.5mm]
&$+e^{u}[u_{x}\cos(pu)+u_{y}\sin(pu)]+\sigma e^{u}$&$$& $$
\\ [0.5mm]
\hline 13&
$u_{t}=(u^{k}u_{x})_{x}+(u^{k}u_{y})_{y}+$&$<\partial_{t}, \partial_{x}, \partial_{y}, (m-k)D_{0}+$&
$m\neq k, $\\
&$+u^{m}[u_{x}\cos(p\ln u)+$&$+D_{3}+pJ_{12}>$&$(m,p)\neq (0,0)$
\\ [0.5mm]
&$+u_{y}\sin(p\ln u)]+\sigma u^{2m-k+1}$&$$&$$
\\ [0.5mm]
\hline 14& $u_{t}=(u^{k}u_{x})_{x}+(u^{k}u_{y})_{y}+$&$<\partial_{t}, \partial_{x}, \partial_{y}, $&$k\neq0$\\
&$+\ln uu_{x}+\sigma u^{-k+1}$&$kD_{0}-D_{3}-t\partial_{x}>$&$$
\\ [0.5mm]
\hline 15&
$u_{t}=(u^{k}u_{x})_{x}+(u^{k}u_{y})_{y}+$ & $<\partial_{t}, \partial_{x}, \partial_{y},D_{3}+pJ_{12}>
$& $k\neq0, p\neq0$\\ [0.5mm]
&$+u^{k}[u_{x}\cos(p\ln u)+$&&
\\ [0.5mm]
&$+u_{y}\sin(p\ln u)]+\sigma u^{k+1}$&$$&$$
\\ [0.5mm]
\hline 16& $u_{t}=(u^{k}u_{x})_{x}+(u^{k}u_{y})_{y}+$ & $<\partial_{t}, \partial_{x}, \partial_{y}, D_{3}>$&
$k\neq0$,\\ [0.5mm]
&$+u^{k}u_{x}+\lambda_{3} u^{k+1}$&&$k^{2}\lambda_{3}\neq 4(k+1)$
\\ [0.5mm]
\hline 17& $u_{t}=(u^{k}u_{x})_{x}+(u^{k}u_{y})_{y}+$ & $<\partial_{t}, \partial_{x}, \partial_{y}, D_{3}, R_{1}, R_{2}>$& $k\neq -1; 0$
\\ [1mm]
&$+4\frac{k+1}{k}u^{k}u_{x}+4\frac{k+1}{k^{2}}u^{k+1}$&$$& $$
\\ [1mm]
\hline 18 & $u_{t}=\triangle u+uu_{x}+\gamma_{1}u$ & $<\partial_{t}, \partial_{x}, \partial_{y}, {\cal G}_{1}>$ & $$ \\ [0.5mm]
\hline 19 & $u_{t}=\triangle u+uu_{x}$ & $<\partial_{t},
\partial_{x}, \partial_{y}, D_{0}-u\partial_{u}, G_{0}>$& $$\\ [0.5mm]
\hline 20& $u_{t}=\triangle u+\ln uu_{x}$ & $<\partial_{t}, \partial_{x}, \partial_{y}, G_{1}>$&
$$  \\ [0.5mm]
\hline 21 & $u_{t}=\triangle u+\ln uu_{x}+\gamma_{1}u\ln u$ & $<\partial_{t}, \partial_{x}, \partial_{y},  {\cal G}_{2}>$&$$
\\ [0.5mm]
\hline 22 & $u_{t}=\triangle u+2\gamma_{1}\ln uu_{x}+$ &
$<\partial_{t},
\partial_{x}, \partial_{y},Y>$& \\ [0.5mm]
&$+u(\ln^{2} u+q)$&&
\\ [0.5mm]
\hline
\end{longtable}

\textbf{Sketch of the proof.}

If one compare the equations and MAIs listed in  Table  \ref{Tabl2-2} with those from  Table  \ref{Tabl4-1}
then can be identified that the  20  cases are identical in the  both tables. In fact,
Cases 1--3,  6--8, 11, 12, 15, 17--19, 21--24, 26,  28, and 31--32
from   Table  \ref{Tabl2-2}
are exactly  Cases  1--3, 4--6,  7, 8, 9, 10--12,  13--16, 17, 18, and   21--22
from  Table  \ref{Tabl4-1}, respectively.


The remaining 12 cases from Table  \ref{Tabl2-2} (Cases  4--5, 9--10, 13--14, 16, 20, 25,  27  and 29--30) are reducible to 11 cases in Table  \ref{Tabl4-1} using the relevant FPTs as  listed in Table  \ref{Tabl5-1}. As one may note making a simple analysis of  Table  \ref{Tabl5-1}, only two additional cases (see the last two lines) should be added to Table  \ref{Tabl4-1} (see 19 and 20 therein). As a result Table  \ref{Tabl4-1} contains exactly 22 equations.

 The sketch of the proof is now completed.\quad $\Box$

 \begin{remark}  The (1+2)-dimensional Burgers equation
 \[  u_{t}=\triangle u+uu_{x} +  uu_{y}, \]
 which is a natural two-dimensional generalization of the famous Burgers equation  $u_{t}=u_{xx}+uu_{x}$,
 is obtainable from the equation listed in Case 19 of Table \ref{Tabl4-1} by the equivalence transformation
 \[t\rightarrow t,\; x\rightarrow \frac{1}{2}(x+y),\; y\rightarrow \frac{1}{2}(x-y),\; u\rightarrow u.\]
 Lie symmetries of the (1+2)-dimensional Burgers equation (in the form listed in Case 19 of Table \ref{Tabl4-1}) were found for the first time in \cite{ed-broa-94} while symmetry reductions and exact solutions are presented in \cite{ed-broa-95}.
 \end{remark}

\newpage

\begin{longtable}[�]{|c|c|c|c|}
\caption[�]{ Mapping  the nonlinear  RDC equations from  Table  \ref{Tabl2-2} to their canonical forms using FPTs from Table  \ref{Tabl3-1}.  The numbers in three columns refer to the relevant cases in Tables  \ref{Tabl2-2}, \ref{Tabl3-1} and \ref{Tabl4-1} }
\label{Tabl5-1}\\
\hline  RDC equation   &  FPTs & Canonical form
\\ [0.5mm]
  in Table  \ref{Tabl2-2} &in Table  \ref{Tabl3-1}  & in Table  \ref{Tabl4-1} \\ [0.5mm]
\hline  4&1& 3
\\ [0.5mm]
\hline  5&2& 3
\\ [0.5mm]
\hline  9& 3& 5
\\ [0.5mm]
\hline  10&3& 6
\\ [0.5mm]
\hline 13&
4&7
\\ [0.5mm]
\hline  14&4& 8(with  $m=k+1$)
\\ [0.5mm]
\hline 16&
4(with  $k=-1$)& 9
\\ [0.5mm]
\hline  20&
5& 10(with  $p=0, \ m=\delta=1$)
\\ [0.5mm]
\hline  25&
6& 16
\\ [0.5mm]
\hline  27&
7& 17
\\ [0.5mm]
\hline  29&
8& 19
\\ [0.5mm]
\hline  30&
9& 20
\\ [0.5mm]
\hline
\end{longtable}

\section {\bf Examples of exact solutions of a generalization of the porous-Fisher equation }

Here we examine the nonlinear equation
\begin{equation}\begin{array}{l}\label{t0}
u_{t}=(uu_{x})_{x}+(uu_{y})_{y}+\lambda_1u u_{x}+ \lambda u(1-u)
\end{array}\end{equation}
where  $\lambda_1$  and $\lambda$  are arbitrary constants.
Eq. (\ref{t0})  with $\lambda_1=0$  coincides with the so called porous-Fisher equation (see, e.g.,  \cite{newman-83, witelski-95, mu-02, McCue_et_al-19}  and its generalization on reaction-diffusion systems \cite{ch-king4} ), which is a generalization of the famous 2D Fisher equation \cite{fi-37}
\begin{equation}\label{t00}
u_{t}= u_{xx}+u_{yy}+ \lambda u(1-u),  \quad \lambda>0.
\end{equation}

Physically Eq. (\ref{t0}) with $\lambda_1=0$ describes  the population dispersing to regions of lower density
more rapidly as the population gets more crowded and has been extensively studied in the 1D approximation. (see, e.g., \cite{mu-02,ch-du-01}
and references therein).

On the other hand,  Eq. (\ref{t0}) can be thought as a generalization of the Murray equation
\begin{equation}\begin{array}{l}\label{t01}
u_{t}= u_{xx}+u_{yy} +\lambda_1u u_{x}+ \lambda u(1-u)
\end{array}\end{equation}
which was intensively studied in \cite{mu-02,ch-se-pl-2018, ch-07} in the 1D approximation.

 It can be noted that Eq. (\ref{t0}) with  $\lambda_{1}=-\lambda=8$ is nothing else but the RDC equation listed in Case 27 of Table \ref{Tabl2-2} under the restriction $k=1$. Moreover, this equation can be simplified via  FPT listed in Case 7 of Table \ref{Tabl3-1}, namely
\begin{equation}\begin{array}{l}\label{t001}
t\rightarrow-\frac{e^{-8t}}{8}, \  x\rightarrow x, \ y\rightarrow y, \  u\rightarrow e^{-8t}u \end{array}\end{equation}
to the form
\begin{equation}\label{t1}
u_{t}=(uu_{x})_{x}+(uu_{y})_{y}+8u u_{x}+ 8u^{2}.
\end{equation}

Now we realize that Eq. (\ref{t1}) possesses the six-dimensional Lie algebra of invariance $AL_6$ (see Case 17 in Table \ref{Tabl4-1}), which is the largest  for   nonlinear RDC equations with non-zero convection terms.
This Lie algebra   is generated by the basic operators
\begin{equation}\begin{array}{l}\label{t2}
\medskip
\partial_{t},\, \partial_{1},\, \partial_{2},\,
D_{3}=t\partial_{t}-u\partial_{u},\\
R_{1}=e^{-x}(\cos y\partial_{x}-\sin y\partial_{y}
-\frac{2}{k}\cos y
u\partial_{u}),\,
R_{2}=e^{-x}(\sin y\partial_{x}+\cos y\partial_{y}
-\frac{2}{k}\sin y
u\partial_{u}).
\end{array}\end{equation}
 It should be noted that Eq.(\ref{t1})  admits  two operators, $R_1$ and $R_2$, with very unusual structure
 and there are not (1+1)-dimensional RDC equations  admitting such kind of operators.

It is well-known that Lie symmetries allow to reduce the given PDE to that of lower dimensionality.
In the case of Eq.(\ref{t1}), there are a wide range of possibilities to make such reductions because
the equation in question admits the six-dimensional Lie algebra of invariance.
Generally speaking, one should construct the so called optimal  systems of inequivalent (non-conjugate) subalgebras of $AL_6$
(see for details \cite{pa-win-77, pa-win-za-75, la-spi-sto-02}).
It is a nontrivial problem in the case of Lie algebras of high dimensionality  and its solving lies beyond scopes of this paper.

On the other hand, there is a straightforward technique  for deriving a set of  reduced equations using the known Lie symmetry of the given PDE.
In the case of Eq.(\ref{t1}), one should take the most general form of Lie's operator belonging to $AL_6$:

\begin{equation}\begin{array}{l}\label{t4}
X=d_{0}\partial_{t}+d_{1}\partial_{x}+d_{2}\partial_{y}+
c_{0}D_{3}+c_{1}R_{1}+c_{2}R_{2}
\end{array}\end{equation}
(here coefficients are arbitrary parameters) and solve the corresponding invariance surface condition
\begin{equation}\begin{array}{l}\label{t420}
X\Phi|_{\Phi=0}=0,
\end{array}\end{equation}
where $\Phi= u(t,x,y)-u_0(t,x,y)$  and  $u_0(t,x,y)$ is an arbitrary solution of Eq. (\ref{t1}). As a result, we arrive at a linear first-order PDE, which is equivalent to the system of three ODEs
\begin{equation}  \label{t4*}\begin{array}{l}
\frac{dt}{c_{0}t+d_{0}}=
\frac{dx}{e^{-x}z(y)+d_{1}}=
\frac{dy}{e^{-x}\dot{z}(y)+d_{2}}
=\frac{du}{-(\frac{2}{k}e^{-x}z(y)+c_{0})u},
\end{array}\end{equation}
where the notation $z(y)=c_{1}\cos y+c_{2}\sin y$ is introduced.
Obviously that the form of solutions of system (\ref{t4*}) depends essentially on the six parameters arising therein.

 Here we examine the $c_{0}=d_{0}=d_{1}=d_{2}=0$ corresponding the symmetry reduction via the operators $R_1$ and $R_2$.
  In this case, the direct integration of system (\ref{t4*}) leads to the first integrals
\begin{equation}  \label{t261}\begin{array}{l}
J_{1}=t,\; J_{2}=\dot{z}(y)e^{x},\; J_{3}=e^{2x}u.
\end{array}\end{equation}
As a result, we arrive at the ansatz
\begin{equation}\label{t5*}
u=e^{-2x}\varphi(t,\omega),\;
\omega=\dot{z}(y)e^{x},
\end{equation}
where  $\varphi=\varphi(t,\omega)$ is a new unknown function.
Substituting  ansatz (\ref{t5*}) into Eq.(\ref{t1}), one obtains the reduced equation
\begin{equation}\begin{array}{l}\label{t5}
\varphi_{t}=(c_{1}^{2}+c_{2}^{2})(\varphi\varphi_{\omega})_{\omega},
\end{array}\end{equation}
which is reducible to the form
\begin{equation}\begin{array}{l}\label{t6}
\varphi_{t}=(\varphi\varphi_{\omega})_{\omega}.
\end{array}\end{equation}
by the time-scaling $t \to \frac{t}{c_{1}^{2}+c_{2}^{2} }. $
 We set $c_{1}^{2}+c_{2}^{2}=1$ in what follows, hence $z(y)=\sin(y+y_0), \ y_0 \in \mathbb{R}$.
Thus, the symmetry reduction of Eq.(\ref{t1}) via the operators $R_1$ and $R_2$ gives nothing else but the porous diffusion equation (\ref{t6}), which is often  called the Boussinesq equation.

Because the Boussinesq equation was extensively studied by many authors (see e.g. Section 4.2 in \cite{ch-se-pl-2018} and references therein) its  exact solutions  were constructed in several works and practically all of them are summarized in the handbook \cite{po-za-04}.

For example, the Boussinesq equation possesses the  plane wave solution
\begin{equation} \label{t11}
\varphi =
p(\omega +pt) +c_3
\end{equation}
and
\begin{equation} \label{t14}
\varphi =c_4t^{-1/3}  - \frac{1}{6}
t^{-1} \omega^2,
\end{equation}
which are obtainable via the  further Lie symmetry reduction of Eq. (\ref{t6}) to ODEs (this equation  admits four-dimensional MAI \cite{ov-80}). Notably the exact solution \eqref{t14} for the first time was derived in \cite{baren-52}.

Thus, using ansatz (\ref{t5*}) and solutions \eqref{t11} and  \eqref{t14}, we obtain exact solutions
 \begin{equation}\label{t23}\begin{array}{l}
u=e^{-2x-8t} \left[p
\cos(y+y_0)e^{x}
-\frac{p^{2}}{8}e^{-8t}+
c_{3}\right];
\end{array}\end{equation}
and
\begin{equation}\label{t28}
u=\frac{4}{3}\left[[
\cos^2(y+y_0)+c_{5}e^{-\frac{16}{3}t-2x}\right]
\end{equation}
of the  nonlinear RDC equation \eqref{t0} with  $\lambda_{1}=-\lambda=8$ (here $c_{5}=-\frac{3}{2}c_{4}(c_{1}^{2}+c_{2}^{2})^{-\frac{1}{3}}$ can be thought as new arbitrary constant).
It should be stressed that we have constructed two  nontrivial exact solutions of \eqref{t0} using the relatively simple solutions \eqref{t11} and  \eqref{t14} of the the Boussinesq equation.

The both exact solutions have interesting asymptotical  behaviour. In fact,  solution  \eqref{t11} tends to zero provided $ t \to \infty$. It means that the solution describes extinction of particles (population of spices, cells  etc.). Solution  \eqref{t11} possesses another time asymptotic because $u \to \frac{4}{3}[
\cos^2(y+y_0)]$ as $ t \to \infty$. So, the solution describes such processes, which tend to the periodical steady-state w.r.t. the variable $y$. Notably both solutions are periodic w.r.t. the variable $y$ and have exponential growth/decay w.r.t. the variable $x$.

Finally, it should be noted that exact solutions of Eq.  \eqref{t0} in the 1D approximation  were constructed and analyzed for the first time in \cite{ch-pl-07} (see  \cite{ch-se-pl-2018} for more details).

 \section {\bf Concluding remarks}

{\it  The main result of this work  consists of solving the LSC problem for the  class of RDC equations of the form (\ref{l1}). The class    contains as particular cases  several subclasses of RDC equations, which has been examined  in earlier papers \cite{go-27, ni-72, na-70, do-kn-sv-83,ed-broa-94,de-iv-soph-08}.  Here this problem is solved for the first time for the most general  class  of such equations (see Theorem 7).  All the specific equations and Lie  symmetries identified  in the above cited  papers follow as particular cases from the results derived in this paper.
 Another important result from applicability point of view consists in deriving
  a new list of point transformations presented in Table \ref{Tabl3-1}. In fact, the  form-preserving transformations listed therein allow us to identified  hidden  relations between  RDC equations. For example, the standard Burgers equation is equivalent to that with a constant source  (see,  case  8 in Table \ref{Tabl3-1}).}

 The LSC problem for  the  class of RDC equations  (\ref{l1}) was solved using two methods (algorithms). The first one is based on the group of ETs and  was firstly  developed and applied  by Ovsiannikov  for (\ref{1-1-1}) with $n=1$ \cite{ov-80}. In the case of class (\ref{l1}), this method lead to 32
  different RDC equations possessing different MIAs  presented in Table \ref{Tabl2-2}
   (notable the first two equations in the table are subclasses of class  (\ref{l1}), however, we treat them as equations containing arbitrary functions as parameters).

 The second method of LSC is based on the set of form-preserving transformations, which do not form any group. Although  FPTs  were firstly  used for solving LSC problem
 in 1990s, the method can be still thought as relatively new (see  Section 2.3 in \cite{ch-se-pl-2018} for more details and references). This method allowed us  to make a further reduction of number of RDC equations possessing nontrivial Lie symmetries. As a result, Theorem \ref{T2-5-6}   was proved, which says that there are  22 RDC equations admitting four- and higher-dimensional MIAs listed in Table \ref{Tabl4-1}.
 It means that there are exactly 22 RDC equations with nontrivial Lie symmetries, which are inequivalent up to any point transformation, and each other RDC  equation possessing a nontrivial Lie symmetry is reducible to one of those from Table \ref{Tabl4-1}. Moreover this list of the RDC equations  cannot be reduced by any point transformations.In this sense these 22  equations form a  canonical list of the RDC  equations with  nontrivial Lie symmetries.

 It should be noted that the results of  LSC   were  verified using the computer algebra package Maple  18.  It means that  each RDC equation  listed in Tables \ref{Tabl2-2} and \ref{Tabl4-1}  was tested and   the  Lie symmetries  obtained coincide with those from Tables   \ref{Tabl2-2} and \ref{Tabl4-1}.

 Now we present a comparison of our results with those derived earlier.
 One notes that the RD equations and the relevant MIAs  listed in cases 2--9 of Table \ref{Tabl4-1} were identified earlier in papers \cite{na-70, do-kn-sv-83}. In particular, it was established for the first time  in \cite{na-70} that the so called conformal exponent $-4/(n+2)$ in the case of two space variables, i.e. $n=2$, leads to infinite-dimensional MIA (see case 9 of Table \ref{Tabl4-1}).

 The diffusion-convection  equations 1--5  listed in Table 1 \cite{ed-broa-94} follow as particular cases from the equations listed in cases 16,15,11,10 and 19 of Table \ref{Tabl4-1}, respectively.

 It should be also stressed that four RDC equations  arising in Cases 10, 12--13  and 15 of Table \ref{Tabl4-1} are absolutely new equations, which have no analogs among (1+1)-dimensional RDC equations (see \cite {ch-se-2006}  or  Table 2.7 in \cite{ch-se-pl-2018} for comparison). The nonlinear  equations  arising in Cases 17--18,   and 20--22 of Table \ref{Tabl4-1}  are also new  (1+2)-dimensional RDC equations. On the other hand, they have analogs among (1+1)-dimensional RDC equations.

 It is quite interesting that Eq.(\ref{t1}) (it is a particular case of  the equation listed in Case 17 of Table \ref{Tabl4-1})  admits six-dimensional Lie algebra $AL_6$ involving two operators, $R_1$ and $R_2$, with very unusual structure. However, calculating Lie brackets, one may show that the subalgebra with the basic operators $\partial_{x}, \partial_{y}, R_1$ and $R_2$ is nothing else but the extended Euclid algebra $AE_1(2)$. In fact, the operators $R_1$ and $R_2$  should be treated as  the space translation operators while $\partial_{y}$ and  $\partial_{x}$ are the operators of rotations and scale transformations, respectively. So, $AL_6$ is the direct sum of $AE_1(2)$  and the two-dimensional Abelian  algebra $\langle \partial_{x}, D_0 \rangle$. It is worth to note that the Lie symmetries $R_1$ and $R_2$ do not occur  as those of (1+1)-dimensional RDC equations. Moreover, we foresee that such symmetries are not allowed by any (1+n)-dimensional RDC equation with $n>2$.

 The Lie symmetries obtained in this work can be widely used for finding exact solutions of the  nonlinear RDC equations arising in applications. It can happen that an equation  possessing a nontrivial Lie symmetry  is not explicitly  listed  in Tables \ref{Tabl2-2} and \ref{Tabl4-1}. However, such equation must  be reducible to one of those in these  tables by appropriate ET or/and FPT.
 In Section 6, an example is presented. In fact, we  have identified  how the porous-Fisher equation with the Burgers term (\ref{t0})(with the  correctly-specified coefficients) can be derived from   Eq.(\ref{t1}). Using the highly nontrivial Lie symmetries of Eq.(\ref{t1}),   exact solutions of the  equation in question  have been constructed and their properties identified. It was also shown that this 2D nonlinear PDE  is reducible to   the 1D Boussinesq equation (\ref{t6}) by the Lie symmetry reduction. The authors are going  to present more  interesting (from the  applicability  point of view) examples  in a forthcoming paper.

It should be also noted that there are  RDC type equations involving the so called gradient-dependent diffusivity (see Eq.(\ref{1-1-0}) with $D(|\nabla u|)$), which arise in many mathematical models describing real-world processes (such as non-Newtonian flow, image processing etc.). The recent papers \cite{ch-ki-ko-16} and \cite{re-anco-17} present new results for solving the LSC problem for  some classes of equations with gradient-dependent diffusivity in the case $n>1$. One may note that Lie symmetry properties of nonlinear  equations with the standard diffusivity  and gradient-dependent diffusivity are essentially different.

 Finally, we point out that the LSC problem (the group classification) for the class of  multidimensional ($n>2$) RDC equations (\ref{1-1-0})  is still open. We foresee that solving  the LSC problem  in the case of three and more  space variables will lead to new results.



\begin{thebibliography}{50}

\bibitem {ames-65} Ames~W F  Nonlinear partial differential equations in engineering.  New York-London: Academic Press;1965.
\bibitem {am-72} Ames~W F Nonlinear partial differential equations in engineering. New York-London: Academic Press;1972.
\bibitem {ko-75} Kozdoba~L Methods for solving nonlinear heat conductions problems. Moscow: Nauka; 1975 (in Russian).
\bibitem {fife-79} Fife~P  Mathematical aspects of reacting and diffusing systems. Berlin-Heidelberg-New York: Springer; 1975.
\bibitem {mu-77} Murray~J D  Nonlinear partial differential equations models in biology. Oxford: Clarendon Press; 1977.
\bibitem {aris-75I} Aris~R  The mathematical theory of diffusion and reaction in permeable catalysts: the theory of the steady state. Oxford: Clarendon Press; 1975.
\bibitem {aris-75II} Aris~R  The mathematical theory of diffusion and reaction in permeable catalysts: Vol. 2: questions of uniqueness stability, and transient behavior. Oxford: Clarendon Press; 1975.
\bibitem {br-03} Britton~N F  Essential mathematical biology. London: Springer; 2003.
\bibitem {edel-05} Edelstein-Keshet~L  Mathematical Models in Biology. Philadelphia:  Society for Industrial and Applied Mathematics; 2005.
\bibitem {ku-na-ei-16} Kuang~Y, Nagy~J D and Eikenberry~S E Introduction to mathematical oncology. Boca Raton, FL: CRC Press; 2016.
\bibitem {mu-02} Murray~J D  Mathematical biology I. New York: Springer; 2002.
\bibitem {mu-03} Murray~J D  Mathematical biology II. New York: Springer; 2003.
\bibitem {ok-le-01} Okubo~A and Levin~S A   Diffusion and ecological problems: modern perspectives. New York: Springer; 2001.
\bibitem {wa-15} Waniewski~J  Theoretical foundations for modeling of mebrane transport in medicine and biomedical engineering. Warsaw: Institute of Computer Science PAS; 2015.
\bibitem {ch-da-2017} Cherniha R and Davydovych V   Nonlinear reaction-diffusion systems -- conditional symmetry, exact solutions and their applications in biology. In: Lecture Notes in Mathematics, vol.2196. Cham:Springer; 2017.
\bibitem{Lie-1881} Lie~S  \"{U}ber die Integration durch bestimmte Integrale von einer Klasse linear partieller Differentialgleichungen. Arch. Math 1881;6(3): 328--68 (in German)
\bibitem{Lie-1885} Lie~S  Algemeine Untersuchungen \"{u}ber Differentialgleichungen, die eine continuirliche endliche Gruppe gestatten. Math. Annalen 1885;25 (in German)
\bibitem{ch-se-98}  Cherniha~R M and Serov~M I  Symmetries, ans\"atze  and exact solutions of  nonlinear second-order evolution equations with convection term. Euro. J. Appl. Math  1998;9: 527--42.
\bibitem{ch-se-2006}  Cherniha~R M and Serov~M I
Symmetries, ans\"atze  and exact solutions of  nonlinear
second-order evolution equations with convection term II. Euro. J.
Appl. Math  2006;17:597--605.
\bibitem {ch-se-pl-2018} Cherniha~R, Serov~M and Pliukhin~O  Nonlinear Reaction-Diffusion-
Convection Equations: Lie and Conditional Symmetry, Exact Solutions and Their Applications. Boca Raton, FL: CRC Press; 2018.
\bibitem{go-27} Goff~J A  Transformations leaving invariant the heart equation of physics. Amer. J. Math 1927;49(1):117-22.
\bibitem{ni-72} Niederer~U  The maximal kinematical invariance group of the free Schr\"{o}dinger equation. Helv. Phys. Acta 1972/73;45(5):802-10.
\bibitem{na-70} Nariboli~G  Self-similar solutions of some nonlinear equations. Appl. Scientific Res  1970;22: 449-61.
\bibitem{do-kn-sv-83} Dorodnitsyn~V A, Knyazeva~I V and Svirshchevskii~S R  Group properties of the nonlinear heart equation with source in the two- and thee-dimensional cases. Differential'niye Uravneniya. 1983;19:1215-23 (in Russian).
\bibitem{ed-broa-94} Edwards~M P and Broadbridge~P   Exact transient solutions to nonlinear diffusion-convection equations in higher dimensions. J. Phys. A: Math. Gen 1994;27: 5455--65.
  \bibitem{de-iv-soph-08}  Demetriou E, Ivanova N M and Sophocleous C  Group analysis of (2+1)- and (3+1)-dimensional diffusion-convection equations. J. Math. Anal. Appl 2008;348:55�65
\bibitem{ol-86}  Olver P. Applications of Lie groups
to differential equations. 2nd ed. In: Graduate Texts in
Mathematics. New York: Springer; 1993.
\bibitem{fu-sh-se-93} Fushchych~W I, Shtelen~W M and Serov~M I  Symmetry analysis and exact solutions of equations of nonlinear mathematical physics. Dordrecht: Kluwer Academic Publishers Group; 1993.
\bibitem{akh91} Akhatov~I S, Gazizov~R K and Ibragimov~N H  Nonlocal symmetries. Heuristic approach. J. Sov. Math  1991;55: 1401-50.
\bibitem{blu-ch-an} Bluman~G W, Cheviakov~A F and Anco~S C  Applications of symmetry
methods to partial differential equations. New York: Springer; 2010.
\bibitem{ga-wi-92} Gazeau~J-P and Winternitz P  Symmetries of variable coefficient Korteweg-de Vries equations. J. Math. Phys  1992;33(12):4087-102.
\bibitem{fi-37}  Fisher~R A  The wave of advance of advantageous genes. Ann. Eugenics 1937;7:353--69.

 \bibitem {newman-83}    Newman~W.I. The long-time behavior of the solution to a non-linear diffusion problem in population genetics and combustion, J. Theor. Biol.1983, 104: 473-484.
     \bibitem {witelski-95} Witelski~T.P.Merging traveling waves for the porous-Fisher�s equation. Appl. Math. Lett.1995,8(4):57�62.
      \bibitem {McCue_et_al-19}   McCue~S.W., Jin Wang, Moroney~T.J. et al. Hole-closing model reveals exponents for nonlinear degenerate diffusivity functions in cell biology. Phys. D 2019; 398: 130�140.
 \bibitem {ch-king4}  Cherniha~R.,  King~J.~R.
  Nonlinear reaction-diffusion systems  with variable diffusivities:
 Lie symmetries, ans\"atze and exact solutions.
   J.  Math. Anal. Appl.  2005, 308:11-35.


    \bibitem {ch-du-01} Cherniha~R  and Dutka V. Exact and Numerical Solutions of the Generalized Fisher Equation.  Rept. Math. Phys 2001;47:393-411.
\bibitem {ch-07} Cherniha~R  New Q-conditional symmetries end exact solutions reaction-diffusion-convection eduations arising in mathematical biology. J. Math. Anal. Appl 2007;326:783-99.
    \bibitem{ed-broa-95} Edwards~M P and Broadbridge~P   Exceptional symmetry reductions of Burgers' equation in
two and three spatial dimensions.Z Angew Math Phys (ZAMP) 1995;46: 595--622.
\bibitem{pa-win-77} Patera~J and Winternitz~P  Subalgebras of real three- and four- dimensional Lie algebras. J. Math. Phys 1977;18(7): 1449-55.
\bibitem{pa-win-za-75} Patera~J, Winternitz~P and Zassenhaus~H  Continuous
subgroups of the fundamental groups of physics I. General method and the Poincare group. J. Math. Phys  1975;16:1597-614.
\bibitem{la-spi-sto-02} Lahno~V I, Spichak~S V and Stognii~V I  Symmetry
analysis of evolutions type equations. Kyiv: Institute of Mathematics of NAS of Ukraine; 2002 (in Ukrainian).
\bibitem{po-za-04} Polyanin~A D and Zaitsev~V F  Handbook of nonlinear partial differential equations. Boca Raton, FL: CRC Press; 2004.
\bibitem{ov-80}   Ovsiannikov LV. The group analysis of differential equations. New York: Academic Press; 1982.
\bibitem{baren-52} Barenblatt~G I  On some unsteady motions of a liquid and gas in a porous medium. Akad. Nauk SSSR. Prikl. Mat. Mekh 1952;16:67--78 (in Russian).
    \bibitem{ch-pl-07} Cherniha R,  Pliukhin O
New conditional symmetries and exact solutions of nonlinear
reaction--diffusion--convection equations.  J. Phys. A: Math. Theor
2007;40:10049--10070.
  \bibitem {ch-ki-ko-16} Cherniha~R,  King~JR  and Kovalenko~S   Lie symmetry properties of nonlinear
reaction-diffusion  equations with gradient-dependent diffusivity.
Commun  Nonlinear Sci. Numer. Simul 2016;36:98-108.
\bibitem {re-anco-17} Recio~E and  Anco~ S  Conservation laws and symmetries of radial generalized nonlinear p-Laplacian evolution equations.  J. Math. Anal. Appl 2017;452:1229�1261.

\end{thebibliography}
\end{document}